\documentclass{amsart}

\usepackage{amsthm}
\usepackage{amsmath}
\usepackage{amssymb}
\usepackage{amsfonts}
\usepackage{latexsym}
\usepackage[all]{xy} \SelectTips{eu}{}
\usepackage{hyperref}
\usepackage{eucal}
\usepackage{xcolor}
\newtheorem{intthm}{Theorem}[]

\newtheorem*{intque*}{Question}
\newtheorem*{intexa*}{Example}
\newcommand{\numberseries}{\bfseries}   
\newlength{\thmtopspace}                
\newlength{\thmbotspace}                
\newlength{\thmheadspace}               
\newlength{\thmindent}                  
\newtheoremstyle{bfupright head,slanted body}
{\thmtopspace}{\thmbotspace}
{\slshape}{\thmindent}{\bfseries}{.}{\thmheadspace}
{{\numberseries \thmnumber{#2\;}}\thmnote{#3}}

\newtheoremstyle{bfupright head,upright body}
{\thmtopspace}{\thmbotspace}
{\upshape}{\thmindent}{\bfseries}{.}{\thmheadspace}
{{\numberseries \thmnumber{#2\;}}\thmnote{#3}}

\newtheoremstyle{fixed bf head,slanted body}
{\thmtopspace}{\thmbotspace}{\slshape}
{\thmindent}{\bfseries}{.}{\thmheadspace}
{{\numberseries \thmnumber{#2\;}}\thmname{#1}\thmnote{ (#3)}}

\newtheoremstyle{fixed bf head,upright body}
{\thmtopspace}{\thmbotspace}{\upshape}
{\thmindent}{\bfseries}{.}{\thmheadspace}
{{\numberseries \thmnumber{#2\;}}\thmname{#1}\thmnote{ (#3)}}

\newtheoremstyle{numbered paragraph}
{\thmtopspace}{\thmbotspace}{\upshape}
{\thmindent}{\upshape}{}{\thmheadspace}
{{\numberseries \thmnumber{#2.}}}

\theoremstyle{bfupright head,slanted body}
\newtheorem{res}{}[section]             \newtheorem*{res*}{}

\theoremstyle{bfupright head,upright body}
\newtheorem{bfhpg}[res]{}               \newtheorem*{bfhpg*}{}

\theoremstyle{fixed bf head,slanted body}
\newtheorem{thm}[res]{Theorem}          \newtheorem*{thm*}{Theorem}
\newtheorem{prp}[res]{Proposition}      \newtheorem*{prp*}{Proposition}
\newtheorem{cor}[res]{Corollary}        \newtheorem*{cor*}{Corollary}
\newtheorem{lem}[res]{Lemma}            \newtheorem*{lem*}{Lemma}
         \newtheorem*{que*}{Question}
\theoremstyle{fixed bf head,upright body}
\newtheorem{dfn}[res]{Definition}       \newtheorem*{dfn*}{Definition}
\newtheorem{rmk}[res]{Remark}           \newtheorem*{rmk*}{Remark}
           \newtheorem*{exa*}{Example}
\newtheorem{setup}[res]{Setup}           \newtheorem*{setup*}{Setup}
\newtheorem{nota}[res]{Notation}           \newtheorem*{nota*}{Notation}

\theoremstyle{numbered paragraph}
\newtheorem{ipg}[res]{}

\newlength{\thmlistleft}        
\newlength{\thmlistright}       
\newlength{\thmlistpartopsep}   
\newlength{\thmlisttopsep}      
\newlength{\thmlistparsep}      
\newlength{\thmlistitemsep}     

\newcounter{eqc}
\newenvironment{eqc}{\begin{list}{\upshape (\textit{\roman{eqc}})}%
		{\usecounter{eqc}%
			\setlength{\leftmargin}{\thmlistleft}%
			\setlength{\labelwidth}{\thmlistleft}%
			\setlength{\rightmargin}{\thmlistright}%
			\setlength{\partopsep}{\thmlistpartopsep}%
			\setlength{\topsep}{\thmlisttopsep}%
			\setlength{\parsep}{\thmlistparsep}%
			\setlength{\itemsep}{\thmlistitemsep}}}%
	{\end{list}}%

\newcounter{prt}
\newenvironment{prt}{\begin{list}{\upshape (\alph{prt})}%
		{\usecounter{prt}%
			\setlength{\leftmargin}{\thmlistleft}%
			\setlength{\labelwidth}{\thmlistleft}%
			\setlength{\rightmargin}{\thmlistright}%
			\setlength{\partopsep}{\thmlistpartopsep}%
			\setlength{\topsep}{\thmlisttopsep}%
			\setlength{\parsep}{\thmlistparsep}%
			\setlength{\itemsep}{\thmlistitemsep}}}%
	{\end{list}}%

\newcounter{rqm}
\newenvironment{rqm}{\begin{list}{\upshape (\arabic{rqm})}%
		{\usecounter{rqm}%
			\setlength{\leftmargin}{\thmlistleft}%
			\setlength{\labelwidth}{\thmlistleft}%
			\setlength{\rightmargin}{\thmlistright}%
			\setlength{\partopsep}{\thmlistpartopsep}%
			\setlength{\topsep}{\thmlisttopsep}%
			\setlength{\parsep}{\thmlistparsep}%
			\setlength{\itemsep}{\thmlistitemsep}}}%
	{\end{list}}%

	{\end{list}}%

\newenvironment{prf*}[1][Proof]{%
	\begin{proof}[\bf #1]
		\setcounter{equation}{0}
		}
	{\end{proof}
}

\newcommand{\pgref}[1]{\ref{#1}}

\renewcommand{\eqref}[1]{(\pgref{eq:#1})}

\newcommand{\eqclbl}[1]{{\upshape(\textit{#1})}}
\newcommand{\proofofimp}[3][:]{\mbox{\eqclbl{#2}$\!\implies\!$\eqclbl{#3}#1}}
\numberwithin{equation}{res}
\def\urltilda{\kern -.15em\lower .7ex\hbox{\~{}}\kern .04em}

\newcommand{\setof}[3][\mspace{1mu}]{\{#1#2 \mid #3#1\}}

\newcommand{\FC}[1]{\mathsf{FlatCot}(#1)}

\newcommand{\PGF}[1]{\mathsf{PGF}(#1)}

\newcommand{\Flat}[1]{\mathsf{Flat}(#1)}
\newcommand{\Cot}[1]{\mathsf{Cot}(#1)}
\newcommand{\Mod}[1]{{_{#1}\mathsf{Mod}}}
\newcommand{\Prj}[1]{\mathsf{Prj}(#1)}
\newcommand{\Inj}[1]{\mathsf{Inj}(#1)}

\newcommand{\xra}[2][]{\xrightarrow[#1]{\:#2\:}}

\newcommand{\Hom}[3][]{\operatorname{Hom}_{#1}(#2,#3)}

\newcommand{\Ext}[4][]{\operatorname{Ext}_{#1}^{#2}(#3,#4)}
\newcommand{\Tor}[4][]{\operatorname{Tor}_{#2}^{#1}(#3,#4)}
\newcommand{\Coker}[1]{\nobreak{\operatorname{Coker}#1}}
\newcommand{\hH}[2][1]{\mathbb{H}_{#1}^{\textnormal{\tiny[}#2\textnormal{\tiny]}}}
\newcommand{\cH}[2][1]{\mathbb{H}^{#1}_{\textnormal{\tiny[}#2\textnormal{\tiny]}}}
\newcommand{\A}{\mathsf{A}}
\newcommand{\Add}{\mathrm{Add}}
\newcommand{\Fun}{\mathrm{Fun}}

\newcommand{\CC}{\mathsf{C}}
\newcommand{\E}{\mathsf{E}}

\newcommand{\G}{\mathsf{G}}
\newcommand{\Se}{\mathsf{S}}

\newcommand{\WG}{\mathsf{W}_{\mathsf{G}}}
\newcommand{\R}{\mathsf{R}}
\newcommand{\X}{\mathsf{X}}

\newcommand{\W}{\mathsf{W}}
\newcommand{\ZZ}{\mathbb{Z}}

\newcommand{\C}{\mathcal{C}}

\newcommand{\Q}{\mathsf{Q}}
\newcommand{\pd}{\mbox{\rm pd}}

\newcommand{\id}{\mathrm{id}}
\newcommand{\kk}{\Bbbk}
\newcommand{\is}{\cong}

\newcommand{\lMod}[1]{{}_{#1}\mspace{-1mu}{\sf{Mod}}}

\newcommand{\rMod}[1]{{\sf{Mod}}_{#1}}

\newcommand{\Cop}{\mathcal{C}^{\mathsf{op}}}
\newcommand{\Aop}{A^{\mathsf{op}}}
\def\soft#1{\leavevmode\setbox0=\hbox{h}\dimen7=\ht0\advance
	\dimen7 by-1ex\relax\if t#1\relax\rlap{\raise.6\dimen7
		\hbox{\kern.3ex\char'47}}#1\relax\else\if T#1\relax
	\rlap{\raise.5\dimen7\hbox{\kern1.3ex\char'47}}#1\relax
	\else\if d#1\relax\rlap{\raise.5\dimen7\hbox{\kern.9ex
			\char'47}}#1\relax\else\if D#1\relax\rlap{\raise.5\dimen7
		\hbox{\kern1.4ex\char'47}}#1\relax\else\if l#1\relax
	\rlap{\raise.5\dimen7\hbox{\kern.4ex\char'47}}#1\relax
	\else\if L#1\relax\rlap{\raise.5\dimen7\hbox{\kern.7ex
			\char'47}}#1\relax\else\message{accent \string\soft
		\space #1 not defined!}#1\relax\fi\fi\fi\fi\fi\fi}

\begin{document}

\title[Flat models for Q-shaped derived categories via PGF objects]%
{Flat models for Q-shaped derived categories via PGF objects}

\author[Z.X. Di]{Zhenxing Di}

\address{Zhenxing Di: School of Mathematical Sciences, Huaqiao University, Quanzhou 362021, China}

\email{dizhenxing@163.com}

\author[L.P. Li]{Liping Li}

\address{Liping Li: Department of Mathematics, Hunan Normal University, Changsha 410081, China}

\email{lipingli@hunnu.edu.cn}

\author[L. Liang]{Li Liang}

\address{Li Liang: Department of Mathematics, Gansu Center for Fundamental Research in Complex Systems Analysis and Control, Lanzhou Jiaotong University, Lanzhou 730070, China}

\email{lliangnju@gmail.com}

\urladdr{https://sites.google.com/site/lliangnju}

\author[Y.J. Ma]{Yajun Ma}

\address{Yajun Ma: Department of Mathematics, Gansu Center for Fundamental Research in Complex Systems Analysis and Control, Lanzhou Jiaotong University, Lanzhou 730070, China}

\email{yjma@mail.lzjtu.cn}

\thanks{Z.X. Di was partly supported by NSF of China (Grant No. 12471034),
Scientific Research Fund of Fujian Province (Grant No. 605-52525002) and
Scientific Research Fund of Huaqiao University (Grant No. 605-50Y22050); L.P. Li was partly supported by NSF of China (Grant No. 12171146); L. Liang was partly supported by NSF of China (Grant No. 12271230) and NSF of Gansu Province (Grant No. 26JRRA059).}

\keywords{Q-shaped derived category, abelian model structure, projectively coresolved Gorenstein flat object.}

\renewcommand{\thefootnote}{\alph{footnote}}
\footnotetext{2020 \emph{Mathematics Subject Classification}. 18G25; 18A25.}

\begin{abstract}
We develop a unified approach, based on projectively coresolved Gorenstein flat (PGF) objects, for constructing flat model structures on diagram categories. Specifically, we show that PGF objects in such categories are fully determined by their objectwise components, which in turn enables us to establish hereditary abelian model structures whose trivial cofibrant objects are precisely the flat objects. As an application, we reobtain flat model structures on $\Q$-shaped derived categories, thereby providing a common framework that subsumes classical constructions for chain complexes. Moreover, we obtain an explicit description of the cofibrant objects in these models.
\end{abstract}

\maketitle
\vspace*{.3cm}
\vspace*{-.9cm}
\enlargethispage{.4cm}

\section*{Introduction}
\noindent
Model structures on diagram categories provide a powerful framework for constructing and studying homotopy categories beyond the classical realm of chain complexes. In particular, Holm and J{\o}rgensen \cite{HJ21} established projective and injective model structures on categories of $\Q$-shaped diagrams, thereby extending the classical model structures for chain complexes to a substantially broader class of diagram categories. This line of work was later generalized in \cite{SV26}, where the projective/injective model structures given in \cite{HJ21} were revisited from a more general perspective.

Motivated by these developments, in this paper we develop a method for constructing flat model structures on diagram categories via projectively coresolved Gorenstein flat objects, or PGF objects for short, introduced in \cite{SS20} by \v{S}aroch and \v{S}t'ov'{\i}\v{c}ek. More precisely, we investigate the behavior of PGF objects in functor categories and employ them to build flat abelian model structures. As an application, we apply our construction to categories of $\Q$-shaped diagrams and reobtain a flat model structure given in \cite{DLLM}, whose homotopy category coincides with the $\Q$-shaped derived category introduced in \cite{HJ21}.

\textbf{Throughout this paper}, we use the term \emph{``subcategory''} to mean a full additive subcategory closed under isomorphisms and finite coproducts. Let $\kk$ be a commutative ring, $A$ an associative $\kk$-algebra, and $\C$ a small $\kk$-linear category such that $\C(p,q)$ is a projective $\kk$-module for each pair of objects $p,q\in\C$. Denote by $\lMod{\C,A}$ the category of $\kk$-linear functors from $\C$ to the category $\lMod{A}$ of (left) $A$-modules, and by $\rMod{\C,A}$ the category of $\kk$-linear functors from $\C^{\sf op}$ to the $\Aop$-module category $\rMod{A}$.

The following result provides, via PGF objects, a method for constructing flat abelian model structures on $\lMod{\C,A}$, that is, model structures whose trivial cofibrant objects are precisely the flat objects, and also gives a characterization for flat objects in $\lMod{\C,A}$; see Theorems \ref{thm:model structure-PGF} and \ref{flat2}.

\begin{intthm}\label{thmA}
Assume that $\C$ is Hom-finite, locally bounded and has a Serre functor (see \ref{setup}). Then the following statements hold.
\begin{prt}
\item Given a set $\G$ of PGF objects in $\lMod{\C,A}$, there exists a hereditary
      abelian model structure
      \[({^{\perp}(\Cot{\C,A}\cap\WG)}, \, \WG, \, \Cot{\C,A})\]
      on $\lMod{\C,A}$, where the trivial cofibrant objects are precisely the flat objects in $\lMod{\C,A}$.

\item If, in addition, $\C$ satisfies the strong retraction property such that the
      pseudo-radical $\mathfrak{r}$ is nilpotent, then an object $X$ in $\lMod{\C,A}$ is flat if and only if $\hH[1]{q}(X)=0$ and $C_q(X)$ is a flat $A$-module for each $q\in\C$.
\end{prt}
\end{intthm}

In the above result, the homology functor $\hH[1]{q}(-)$ is defined as $\Tor[\C]{1}{S\{q\}}{-}$, and the cokernel functor $C_q(-)$ is defined as $S\{q\} \otimes_\C X$, where $S\{q\}: \Cop\to\lMod{\kk}$ denotes the stalk functor; see \ref{adjoint triple2}. We note that statement (b) in Theorem \ref{thmA} remains valid without the assumptions that $\C$ is locally bounded and has a Serre functor; see Theorem \ref{flat2}. This complements \cite[Theorem E]{HJ21}, where Holm and J{\o}rgensen established analogous results for projective and injective objects.

To apply Theorem \ref{thmA}, we need to identify a set of PGF objects in $\lMod{\C,A}$. The following result shows that PGF objects in $\lMod{\C,A}$ are completely determined by their objectwise PGF components in $\Mod{A}$.

\begin{intthm}\label{thmB}
Assume that $\C$ satisfies the conditions specified in Theorem \ref{thmA}. Then an object $F$ in $\lMod{\C,A}$ is a PGF object if and only if $F(q)$ is a PGF $A$-module for each $q\in\C$.
\end{intthm}

Theorem \ref{thmB} is proved using the tool of left Frobenius pairs, which correspond to the weak AB-contexts introduced in \cite{MAsROB89}. Specifically, we establish a method for constructing left Frobenius pairs in $\lMod{\C,A}$ by showing that $(\Fun(\C,\X), \widetilde{\W})$ is a left Frobenius pair in $\Mod{\C,A}$ whenever $(\X, \W)$ is a left Frobenius pair in $\Mod{A}$; see Theorem \ref{Fpair}. In this way, we obtain a left Frobenius pair $(\Fun(\C,\PGF{A}), \widetilde{\Prj{A}})$ in $\Mod{\C,A}$ (see Corollary \ref{Fpair-exa}), which is then used to prove Theorem \ref{thmB}.

\begin{equation*}
  \ast \ \ \ast \ \ \ast
\end{equation*}

We apply our main results to the category $\lMod{\Q,R}$ of $\Q$-shaped diagrams studied in \cite{HJ21}. By definition, $\Q$ is a small preadditive category that is Hom-finite, locally bounded, has a Serre functor, and satisfies the strong retraction property with pseudo-radical $\mathfrak{r}$ nilpotent. By Theorem \ref{thmB}, $S_q(R)=S\langle q\rangle\otimes_\ZZ R$ is a PGF object in $\lMod{\Q,R}$. Taking $\G=\{S_{q}(R)\mid q\in \Q\}$ and applying Theorem \ref{thmA}, we reconstruct the abelian model structure on $\lMod{\Q,R}$ given in \cite{DLLM}. The following result is a special case of Theorem \ref{flat model structure}.

 \begin{intthm}\label{thmC}
Let $\Q$ be a small preadditive category that is Hom-finite, locally bounded, has a Serre functor and satisfies the strong retraction property. Assume that the pseudo-radical $\mathfrak{r}$ is nilpotent, and $R$ is any ring. Then there exists a hereditary abelian model structure
\[
{\mathsf M}=(^{\perp}(\Cot{\Q,R}\cap\E), \, \E, \, \Cot{\Q,R})
\]
on $\lMod{\Q,R}$, where $\E$ is the subcategory of exact objects in $\lMod{\Q,R}$. Consequently, its homotopy category coincides with the $\Q$-shaped derived category introduced in \cite[Corollary 2.8]{HJ21}.
 \end{intthm}

Our next result gives an explicit description of the cofibrant objects in the above abelian model structure ${\mathsf M}$, as a special case of Proposition \ref{cofibrant} and Theorem \ref{cofibrant object}.

\begin{intthm}\label{thmD}
Let $\Q$ be as in Theorem $\ref{thmC}$ and $R$ any ring. If $X\in{^\perp(\Cot{\Q,R}\cap\E)}$, then $X(q)$ is a flat $R$-module for each $q\in\Q$. Moreover, the converse holds whenever $R$ has finite weak global dimension.
\end{intthm}

We mention that Theorem \ref{thmD} is known when $\Q$ is the path category of
\[
\Gamma=\cdots \to \bullet \xra{\partial} \bullet \xra{\partial} \bullet \to \cdots
\]
subject to the relations $\partial^2=0$. Indeed, in this case, $\lMod{\Q,R}$ may be identified with the category of chain complexes, and ${^\perp(\Cot{\Q,R}\cap\E)}$ is precisely the subcategory of dg-flat complexes; see \cite{Gi04}.

\section{Preliminaries}\label{property}
\noindent
In this section, we mainly collect some basic results on homological objects in $\lMod{\C,A}$ that will be used frequently throughout the paper.

\begin{bfhpg}[\bf Cotorsion pairs]\label{cotpair}
Let $\A$ be an abelian category.
A pair $(\CC, \W)$ of subcategories of $\A$ is called a \emph{cotorsion pair} if $\CC^{\bot}=\W$ and $^{\bot}\W=\CC$, where
\begin{align*}
{\CC}^{\bot} &=
\{N\in {\sf A}\ |\ {\rm Ext}_{\A}^{1}(C,N)=0 {\rm \ for\ all}\ C\in\CC\},\ \text{and}\\
^{\bot}{\W}  &=
\{M\in {\sf A}\ |\ {\rm Ext}_{\A}^{1}(M,W)=0 {\rm \ for\ all}\ W\in\W\}.
\end{align*}
For any subcategory $\CC$ of $\A$, the pair $({^\perp(\CC^\perp)},\CC^\perp)$ is a cotorsion pair, called the cotorsion pair \emph{generated} by $\CC$. Recall from Enochs and Jenda \cite{rha} that a cotorsion pair $(\CC,\W)$ is \emph{complete} if for each object $M$ in $\A$, there exist short exact sequences
\[
0\to M\to W\to C\to 0
\quad \text{and} \quad
0\to W'\to C'\to M\to 0
\]
with $W, W'\in\W$ and $C, C'\in\CC$. A cotorsion pair $(\CC,\W)$ is called \emph{hereditary} if $\CC$ is closed under kernels of epimorphisms and $\W$ is closed under cokernels of monomorphisms.
\end{bfhpg}

\begin{bfhpg}[\bf Abelian model structures]
Let $\A$ be an abelian category. A triple $(\Q,\W,\R)$ of subcategories of $\A$ is called a \emph{Hovey triple} if
 \[
 (\Q,\W\cap \R) \quad \text{and} \quad (\Q\cap\W, \R)
 \]
are complete cotorsion pairs in $\A$ with $\W$ thick, that is, it is closed under summands, and such that in any short exact sequence, if two of the three terms lie in $\W$, then so does the third.

In \cite{Ho02}, Hovey established a one-to-one correspondence between abelian model structures on $\A$ and Hovey triples in $\A$. Specifically, from an abelian model structure one obtains a Hovey triple $(\Q, \W, \R)$, where $\Q$ is the class of cofibrant objects, $\R$ is the class of fibrant objects, and $\W$ is the class of trivial objects. Conversely, a Hovey triple $(\Q, \W, \R)$ gives rise to an abelian model structure, where a cofibration (resp., trivial cofibration) is a monomorphism whose cokernel lies in $\Q$ (resp., $\Q\cap\W$), and a fibration (resp., trivial fibration) is an epimorphism whose kernel lies in $\R$ (resp., $\R\cap\W$). By this result, an abelian model structure on $\A$ can be succinctly represented by a Hovey triple.

A Hovey triple $(\Q, \W, \R)$ is called \emph{hereditary} if the cotorsion pairs $(\Q, \W\cap\R)$ and $(\Q\cap\W, \R)$ are hereditary. In that case, the corresponding homotopy category is triangulated and is triangle equivalent to the stable category $(\Q\cap\R)/(\Q\cap\W\cap \R)$; see \cite{Gil162}.
\end{bfhpg}

\begin{bfhpg}[Injective cogenerators]
Let $\A$ be an abelian category. For two subcategories $\W$ and $\X$ of $\A$, recall from \cite{MAsROB89} that $\W$ is a \emph{cogenerator} for $\X$ if $\W\subseteq\X$ and for each $X\in\X$, there is a short exact sequence $0 \to X \to W \to X' \to 0$ with $W\in\W$ and $X'\in\X$. The subcategory $\W$ is said to be an \emph{injective cogenerator} for $\X$ if $\W$ is a cogenerator for $\X$ and $\Ext[\A]{i}{X}{W}=0$ for all $X\in\X$, $W\in\W$ and $i\geqslant1$. Dually, one may give the definition of \emph{projective generators}.
\end{bfhpg}

\begin{bfhpg}[\bf Left Frobenius pairs]\label{fp-def}
Recall from \cite{BMPS19} that a pair $(\X,\W)$ of subcategories of $\A$ is called a \emph{left Frobenius pair} if $\X$ is closed under extensions, kernels of epimorphisms and direct summands in $\A$, and $\W$ is an injective cogenerator for $\X$ and closed under direct summands in $\A$.
\emph{Right Frobenius pairs} are defined dually.

Frobenius's name is invoked here because if $(\X,\W)$ is a left Frobenius pair and $\W$ is a projective generator for $\X$, then $\X$ is a Frobenius category in the usual sense with projective-injective objects precisely the objects in $\W$; see \cite[Remark 2.6]{BMPS19}. We mention that there is a bijective correspondence between left Frobenius pairs and weak AB-contexts introduced in \cite{MAsROB89}; see \cite[Theorem 5.4]{BMPS19}.
\end{bfhpg}

\begin{bfhpg}[\bf Tensor products in $\lMod{\C,A}$]\label{tensor}
By \cite[Proposition 3.12]{HJ21}, $\lMod{\C,A}$ is a locally finitely presentable Grothendieck category with enough projectives. Following Oberst and R\"{o}hrl \cite{Oberst70}, one may define a tensor product functor
\[-\otimes_{\C,A}-: \rMod{\C,A}\times \lMod{\C,A}\to \lMod{\kk}\]
as follows; it is right exact in each variable.

Fix $X\in \lMod{\C,A}$ and let $G$ be a $\Bbbk$-module. Define an object $\Hom[\Bbbk]{X}{G}$ in $\rMod{\C,A}$ as follows:
 \begin{prt}
\item[$\bullet$] For each $q\in\C^{\sf op}$, set
                 $\Hom[\Bbbk]{X}{G}(q)=\Hom[\Bbbk]{X(q)}{G}\in \rMod{A}.$
\item[$\bullet$] For each morphism $f: p\to q$ in $\C^{\sf op}$, define
\[\Hom[\Bbbk]{X}{G}(f)=\Hom[\Bbbk]{X(f^{\sf op})}{G}:\Hom[\Bbbk]{X(p)}{G}\to \Hom[\Bbbk]{X(q)}{G}.\]
\end{prt}
It is evident that $\Hom[\Bbbk]{X}{-}$ is a functor from $\lMod{\Bbbk}$ to $\rMod{\C,A}$. This functor is left exact and preserves arbitrary products, so it admits a left adjoint functor from $\rMod{\C,A}$ to $\lMod{\Bbbk}$, denoted by $-\otimes_{\C,A} X$, which serves as the tensor product in our setting.
\end{bfhpg}

\begin{rmk}\label{useful}
For each $X\in\lMod{\C,A}$, $Y\in\rMod{\C,A}$ and $G\in\lMod{\Bbbk}$, the following natural isomorphism holds (see \cite[Theorem 1]{St68}):
\begin{align*}\label{1.9.1}\tag{\ref{useful}.1}
\Hom[\C^{\sf op},A^{\sf op}]{Y}{\Hom[\Bbbk]{X}{G}}
& \cong \Hom[\Bbbk]{Y\otimes_{\C,A} X}{G}                \\
& \cong \Hom[\C,A]{X}{\Hom[\kk]{Y}{G}}.
\end{align*}
\end{rmk}

\begin{bfhpg}[\bf Purity]
A short exact sequence $0 \to X \to X' \to X'' \to 0$ in $\lMod{\C,A}$ is called \emph{pure} if for every finitely presented object $T$ in $\lMod{\C,A}$, the sequence
\[
0\to \Hom[\C,A]{T}{K} \to \Hom[\C,A]{T}{M} \to \Hom[\C,A]{T}{N} \to 0
\]
is exact. In this case, $X\to X''\to 0$ is called a \emph{pure epimorphism}, and $0\to X \to X'$ is called a \emph{pure monomorphism}.

By \cite[Theorem 2]{St68}, a short exact sequence
$0 \to X\to X'\to X''\to 0$
in $\lMod{\C,A}$ is pure if and only if for each object $Y$ in $\rMod{\C,A}$, the sequence
\[
0 \to Y\otimes_{\C,A} X \to Y\otimes_{\C,A} X' \to Y\otimes_{\C,A} X'' \to 0
\]
is exact in $\lMod{\Bbbk}$, and if and only if the sequence
\[
0 \to \Hom[\Bbbk]{X''}{I} \to \Hom[\Bbbk]{X'}{I} \to \Hom[\Bbbk]{X}{I} \to 0
\]
splits, where $I$ is an injective cogenerator in $\lMod{\Bbbk}$.
\end{bfhpg}

\begin{bfhpg}[\bf Flat objects]\label{representable}
An object $F$ in $\lMod{\C,A}$ is called \emph{flat} if every epimorphism $M\to F$ is a pure epimorphism. It is not difficult to check that $F$ is flat if and only if it is a filtered colimit of finitely generated projective objects.

Let ${\sf Flat(\C,A)}$ denote the subcategory of flat objects in $\lMod{\C,A}$. By Rump \cite[Proposition 7]{Ru10}, ${\sf Flat(\C,A)}$ is closed under extensions, filtered colimits, pure monomorphisms, and pure epimorphisms.

An object $X$ in $\lMod{\C,A}$ is called \textit{cotorsion} if it lies in $\Flat{\C,A}^\perp$. Denote the subcategory of cotorsion objects by $\Cot{\C,A}$.
By \cite[Theorem A.3]{HJ21}, the pair $(\Flat{\C,A}, \Cot{\C,A})$ is a complete hereditary cotorsion pair generated by a set, that is, there exists a set $\Se\subseteq{\sf \Flat{\C,A}}$ such that $\Se^\perp=\Cot{\C,A}$.
\end{bfhpg}

The following result is established in \cite[Theorem 3]{St68}.

\begin{lem}\label{flat-inj}
Let $X$ be an object in $\lMod{\C,A}$ and $I$ an injective cogenerator in $\lMod{\Bbbk}$. Then the following statements are equivalent.
\begin{eqc}
\item $X$ is flat.
\item $\Hom[\Bbbk]{X}{I}$ is injective in $\rMod{\C,A}$.
\item The functor $-\otimes_{\C,A}X$ is exact.
\item $\Tor[\C,A]{1}{-}{X}=0$.
\end{eqc}
\end{lem}

\begin{cor}\label{cotorsion}
Let $X$ be an object in $\lMod{\C,A}$ and $I$ an injective module in $\lMod{\Bbbk}$. Then $\Hom[\Bbbk]{X}{I}$ is cotorsion in $\rMod{\C,A}$.
\end{cor}
\begin{prf*}
It suffices to prove that $\Ext[\C^{\mathrm{op}},A^{\mathrm{op}}]{1}{F}{\Hom[\Bbbk]{X}{I}}=0$ for every flat object $F$ in $\rMod{\C,A}$.
Take a projective resolution $P^{\bullet}$ of $F$ in $\rMod{\C,A}$.
By Lemma \ref{flat-inj}, $\Hom[\Bbbk]{P^{\bullet}}{I}$ is an injective resolution of $\Hom[\kk]{F}{I}$ in $\lMod{\C,A}$. Hence,
\begin{align*}
\Ext[\C^{\mathrm{op}},A^{\mathrm{op}}]{1}{F}{\Hom[\Bbbk]{X}{I}}
&=\mathrm{H}^{1}\Hom[\C^{\mathrm{op}},A^{\mathrm{op}}]{P^{\bullet}}{\Hom[\Bbbk]{X}{I}}\\
&\cong \mathrm{H}^{1}\Hom[\C,A]{X}{\Hom[\Bbbk]{P^{\bullet}}{I}}                   \\
&\cong \Ext[\C,A]{1}{X}{\Hom[\Bbbk]{F}{I}}                                       \\
&= 0,
\end{align*}
where the first isomorphism follows from (\ref{1.9.1}), and the last equality holds as $\Hom[\kk]{F}{I}$ is injective by Lemma \ref{flat-inj}.
\end{prf*}

\begin{ipg}\label{adjoint triple}
By \cite[Corollary 3.9]{HJ21}, for each object $q \in \C$, there is an adjoint triple $(F_{q},E_{q},G_{q})$ as follows:
  \begin{equation*}
  \xymatrix@C=4pc{
      \lMod{\C,A}
    \ar[r]^-{E_{q}}
    & \lMod{A}
    \ar@/_1.8pc/[l]_-{F_{q}}
    \ar@/^1.8pc/[l]^-{G_{q}}
  }
  \qquad \text{given by} \qquad
  {\setlength\arraycolsep{1.5pt}
   \renewcommand{\arraystretch}{1.2}
  \begin{array}{rcl}
  F_q(M) &=& \C(q,-) \otimes_\kk M \\
  E_q(X) &=& X\mspace{1.5mu}(q) \\
  G_{q}(M) &=& \Hom[\kk]{\C(-, q)}{M}\;.
  \end{array}
  }
  \end{equation*}
\end{ipg}

It is clear that the \textit{evaluation} functor $E_q$ is exact. Since $\C(p,q)$ is a projective $\kk$-module for each pair of objects $p$ and $q$ in $\C$ by the assumption, both $F_q$ and $G_q$ are exact as well. Consequently, by Holm and J{\o}rgensen \cite[Lemma 5.1]{HJ19}, we obtain the following result.

\begin{lem}\label{n-iso}
Let $X$ be an object in $\lMod{\C,A}$ and $q$ an object in $\C$. Then for each $M \in \lMod{A}$ and $i \geqslant 0$, there are isomorphisms
\[
\Ext[\C,A]{i}{F_q(M)}{X} \is \Ext[A]{i}{M}{X(q)}
\,\text{  and  }\,
\Ext[A]{i}{X(q)}{M} \is \Ext[\C,A]{i}{X}{G_q(M)}.
\]
\end{lem}

\begin{nota}
Let $\Prj{\C,A}$ (resp., $\Inj{\C,A}$) denote the subcategory of projective (resp., injective) objects in $\lMod{\C,A}$. When $A=\kk$, we write $\lMod{\C}$ (resp., $\Prj{\C}$, $\Inj{\C}$, $\Flat{\C}$ and $\Cot{\C}$) instead of $\lMod{\C,\kk}$ (resp., $\Prj{\C,\kk}$, $\Inj{\C,\kk}$, $\Flat{\C,\kk}$ and $\Cot{\C,\kk}$), and we write $-\otimes_{\C}-$ instead of $-\otimes_{\C,\kk}-$.
\end{nota}

The following result, taken from \cite[Proposition 1.10]{DLLM}, asserts that the evaluation functor preserves certain homological properties.

\begin{prp}\label{degreewise}
Let $X$ be an object in $\lMod{\C,A}$ and $q$ an object in $\C$. Then the following statements hold.
\begin{prt}
\item If $X\in\Prj{\C,A}$, then $X(q)\in\Prj{A}$.
\item If $X\in\Inj{\C,A}$, then $X(q)\in\Inj{A}$.
\item If $X\in\Flat{\C,A}$, then $X(q)\in\Flat{A}$.
\item If $X\in\Cot{\C,A}$, then $X(q)\in\Cot{A}$.
\end{prt}
\end{prp}

The next result, which is an immediate consequence of Lemma \ref{n-iso} and Proposition \ref{degreewise}, shows that the functor $G_q$ preserves cotorsion objects.

\begin{cor}\label{cot}
Let $N$ be a cotorsion $A$-module and $q$ an object in $\C$. Then $G_{q}(N)$ is  cotorsion in $\lMod{\C,A}$.
\end{cor}

\begin{ipg}
Let $(-)^{\natural}$ denote the forgetful functor from $\lMod{\C,A}$ to $\lMod{\C}$. Then $(-)^{\natural}$ is exact and admits the following adjoint triple (see \cite[Corollary 3.5]{HJ21}):
  \begin{equation*}
  \xymatrix@C=4pc{
    \lMod{\C,A}
    \ar[r]^-{(-)^{\natural}}
    &
    \lMod{\C}
    \ar@/_1.8pc/[l]_-{-\otimes_{\kk}A}
    \ar@/^1.8pc/[l]^-{\Hom[\kk]{A}{-}}
  }
  \end{equation*}
\end{ipg}

The following result can be found in \cite[Lemma 1.14]{DLLM}.

\begin{lem}\label{iso}
Assume that $A$ has finite projective dimension over $\kk$. Let $X$ be an object in $\Mod{\C,A}$. Then the following statements hold.
\begin{prt}
\item For each $G \in \Mod{\C}$ admitting an exact sequence $0 \to G \to F_{-1} \to F_{-2} \to \cdots$ with each $F_i \in \Flat{\C}$, there is an isomorphism
\[
\Ext[\C,A]{1}{G \otimes_\kk A}{X} \is \Ext[\C]{1}{G}{X^{\natural}}.
\]
\item For each $G \in \Mod{\C}$ admitting an exact sequence $\cdots \to I_1 \to I_0 \to G \to 0$ with each $I_i \in \Inj{\C}$, there is an isomorphism
\[
\Ext[\C,A]{1}{X}{\Hom[\kk]{A}{G}} \is \Ext[\C]{1}{X^{\natural}}{G}.
\]
\end{prt}
\end{lem}

\begin{setup}\label{setup}
Throughout the paper, we work under one or more of the following assumptions on $\C$ (see \cite{HJ21}):
\begin{rqm}
\item \emph{Hom-finiteness}: Each $\Bbbk$-module $\C(p,q)$ is finitely generated and
      projective.
\item \emph{Local Boundedness}: For every object $q$ in $\C$, there are only finitely many objects in $\C$ mapping nontrivially into or out of $q$.
\item \emph{Existence of a Serre Functor}: There is a $\Bbbk$-linear
      auto-equivalence $\mathbb{S}: \C \to \C$ such that $\C(p,q)\is\Hom[\Bbbk]{\C(q,\mathbb{S}(p))}{\Bbbk}$ for all $p,q\in\C$.
\item \emph{Strong Retraction Property}: For each object $q$ in $\C$, the unit map
      $\Bbbk\to \C(q,q)$, given by $x\mapsto x\cdot\id_q$, admits a $\Bbbk$-module retraction and there exists a collection $\{\mathfrak{r}_q\}_{q\in\C}$ of complements, that is, $\Bbbk$-modules $\mathfrak{r}_q$ such that
      \[\C(q,q)=(\Bbbk\cdot\id_q)\oplus\mathfrak{r}_q,\]
      compatible with composition in $\C$ in the sense that
    \begin{eqc}
    \item[($\dagger$)] $\mathfrak{r}_q\circ\mathfrak{r}_q\subseteq\mathfrak{r}_q$ for all $q\in\C$;
    \item[($\ddagger$)] $\C(q,p)\circ\C(p,q)\subseteq\mathfrak{r}_q$ for all $p\neq q$ in $\C$.
    \end{eqc}
\end{rqm}
\end{setup}

For each $q\in\C$, consider the adjoint triple $(F'_{q},E'_{q},G'_{q})$
  \begin{equation*}
  \xymatrix@C=4pc{
    \rMod{\C,A}
    \ar[r]^-{E'_{q}}
    &
    \rMod{A}
    \ar@/_1.8pc/[l]_-{F'_{q}}
    \ar@/^1.8pc/[l]^-{G'_{q}}
  }
  \qquad \text{given by} \qquad
  {\setlength\arraycolsep{1.5pt}
   \renewcommand{\arraystretch}{1.2}
  \begin{array}{rcl}
  F'_q(N) &=& \C(-,q) \otimes_\kk N \\
  E'_q(Y) &=& Y\mspace{1.5mu}(q) \\
  G'_{q}(N) &=& \Hom[\kk]{\C(q, -)}{N}\;.
  \end{array}
  }
  \end{equation*}
Then the following result holds whenever $\C$ is Hom-finite and has a Serre functor; see the proof of \cite[Lemma 3.8]{DLLM}.

\begin{lem}\label{Tor formula}
Assume that $\C$ is Hom-finite and has a Serre functor $\mathbb{S}$. Then for each $N\in \rMod{A}$,  $X\in\Mod{\C,A}$, $q\in\C$ and $i\geqslant0$, there is an isomorphism
\[
\Tor[i]{\C,A}{G'_{q}(N)}{X}\is\Tor[i]{A}{N}{X(\mathbb{S}(q))}.
\]
\end{lem}

\begin{prf*}
Let $P^{\bullet}$ be a projective resolution of $X$. Since the functor $E_q$ is exact and preserves projectives (as $G_q$ is exact), it follows that  $P^{\bullet}(\mathbb{S}(q))$ is a projective resolution of $X(\mathbb{S}(q))$. Hence,
\begin{align*}
\Tor[i]{\C,A}{G'_{q}(N)}{X}
&\cong \mathrm{H}_{i}(G'_{q}(N)\otimes_{\C,A}P^{\bullet})     \\
&\cong \mathrm{H}_{i}(N\otimes_{A}P^{\bullet}(\mathbb{S}(q))) \\
&\cong \Tor[i]{A}{N}{X(\mathbb{S}(q))},
\end{align*}
where the second isomorphism holds by \cite[Lemma 3.9]{DLLM}.
\end{prf*}

\begin{lem}\label{Tor vanish}
Assume that $\C$ is Hom-finite, locally bounded, and has a Serre functor $\mathbb{S}$. Let $X$ be an object in $\Mod{\C,A}$ and $i\geqslant0$. Then the following statements are equivalent.
\begin{eqc}
\item $\Tor[i]{A}{I}{X(q)}=0$ for each injective module $I$ in $\rMod{A}$ and object $q\in\C$.
\item $\Tor[i]{\C,A}{E}{X}=0$ for each injective object $E$ in $\rMod{\C,A}$.
\end{eqc}
\end{lem}

\begin{prf*}
\proofofimp{i}{ii} Let $E$ be injective in $\rMod{\C,A}$, and let $I$ be a faithfully injective $\Aop$-module. Then for all objects $q\in\C$, the objects $G'_{q}(I)$ are injective and cogenerate $\rMod{\C,A}$; see \cite[Proposition 3.12]{HJ21}. Hence, $E$ is a direct summand of an object of the form $\prod_{q\in\C}G'_{q}(I)$. Moreover, since $\C$ is locally bounded, $\prod_{q\in\C}G'_{q}(I)\cong \oplus_{q\in\C}G'_{q}(I)$ by \cite[Proposition 3.7]{HJ23}. Thus, it suffices to consider $E=G'_{q}(I)$, for which Lemma \ref{Tor formula} gives
$$\Tor[i]{\C,A}{E}{X}\is\Tor[i]{A}{I}{X(\mathbb{S}(q))}=0.$$

\proofofimp{ii}{i} Let $I$ be injective in $\rMod{A}$ and $q$ an object in $\C$. Then $G'_{\mathbb{S}^{-1}(q)}(I)$ is injective in $\rMod{\C,A}$.
Hence, by Lemma \ref{Tor formula},
\[\Tor[i]{A}{I}{X(q)}\is\Tor[i]{\C,A}{G'_{\mathbb{S}^{-1}(q)}(I)}{X}=0,\]
since $G'_{\mathbb{S}^{-1}(q)}$ preserves injectives.
\end{prf*}

\section{A description for flat objects}
\noindent In this section, we provide a description for flat objects in $\lMod{\C,A}$.

\begin{bfhpg}[\bf Pseudo-radical]\label{radical}
Assume that $\C$ satisfies the strong retraction property. The $\Bbbk$-modules
\[
\mathfrak{r}(p,q)=
\begin{cases}
\mathfrak{r}_q, & \text{if } p=q,\\
\C(p,q), & \text{otherwise};
\end{cases}
\]
form an ideal $\mathfrak{r}$ of $\C$, called the \emph{pseudo-radical}. Indeed, viewing $\C$ as a $\kk$-algebra with a complete set of idempotents, $\mathfrak{r}$ is a two-sided ideal of $\C$. Moreover, the quotient $\C/\mathfrak{r}$ is a direct sum of copies of $\kk$ indexed by the objects of $\C$.
\end{bfhpg}

\begin{ipg}\label{adjoint triple2}
Assume that $\C$ satisfies the strong retraction property. By \cite[Proposition 7.15]{HJ21}, for every object $q$ in $\C$, there is an adjoint triple $(C_{q},S_{q},K_{q})$ as follows:
  \begin{equation*}
  \xymatrix@C=4pc{
    \Mod{A}
    \ar[r]^-{S_{q}}
    &
    \Mod{\C,A}
    \ar@/_1.8pc/[l]_-{C_{q}}
    \ar@/^1.8pc/[l]^-{K_{q}}
  }
  \qquad \text{given by} \qquad
  {\setlength\arraycolsep{1.5pt}
   \renewcommand{\arraystretch}{1.2}
  \begin{array}{rcl}
  C_q(X) &=& S\{q\} \otimes_\C X \\
  S_q(M) &=& S\langle q\rangle\otimes_\kk M \\
  K_{q}(X) &=& \Hom[\C]{S\langle q\rangle}{X}\;.
  \end{array}
  }
  \end{equation*}
Here, $S\langle q\rangle: \C\to \Mod{\kk}$ denotes the \emph{stalk functor},
defined by
\[
S\langle q\rangle=\C(q,-)/\mathfrak{r}(q,-)\in\lMod{\C},
\]
and described as follows (see \cite[Lemma 7.10]{HJ21}):
\begin{itemize}
\item for every object $p\in\C$,
\[
S\langle q\rangle(p)=
\begin{cases}
\kk, & \text{if } p=q,\\
0, & \text{otherwise};
\end{cases}
\]
\item for every morphism $f\in\C$,
\[
S\langle q\rangle(f) =
\begin{cases}
x\cdot \id_\kk, & \text{if } f=x\cdot\id_q+g\in(\kk\cdot\id_q)\oplus\tau_q=\C(q,q),\\
0, & \text{otherwise};
\end{cases}
\]
\end{itemize}
The (contravariant) stalk functor $S\{q\}$ is defined dually.
\end{ipg}

\begin{bfhpg}[\bf (Co)homology for objects]
Assume that $\C$ satisfies the strong retraction property. For $X\in\Mod{\C,A}$ and $q\in\C$, define the $i$-th (co)homology of $X$ at $q$ for each $i\geqslant0$ by
\[\cH[i]{q}(X)=\Ext[\C]{i}{S\langle q\rangle}{X}
\quad \text{and} \quad
\hH[i]{q}(X)=\Tor[\C]{i}{S\{q\}}{X}.\]
Note that $\cH[i]{q}$ and $\hH[i]{q}$ are functors from $\Mod{\C,A}$ to $\Mod{A}$.

We mention that there is an adjoint triple $(C'_{q},S'_{q},K'_{q})$ as follows:
  \begin{equation*}
  \xymatrix@C=4pc{
    \rMod{A}
    \ar[r]^-{S'_{q}}
    &
    \rMod{\C,A}
    \ar@/_1.8pc/[l]_-{C'_{q}}
    \ar@/^1.8pc/[l]^-{K'_{q}}
  }
  \qquad \text{given by} \qquad
  {\setlength\arraycolsep{1.5pt}
   \renewcommand{\arraystretch}{1.2}
  \begin{array}{rcl}
  C'_q(Y) &=& Y \otimes_\C S\langle q\rangle \\
  S'_q(N) &=& S\{q\}\otimes_\kk N \\
  K'_{q}(Y) &=& \Hom[\C^{\sf op}]{S\{q\}}{Y}\;.
  \end{array}
  }
  \end{equation*}
Consequently, for each $Y\in\rMod{\C,A}$ and $i\geqslant0$, the $i$-th (co)homology of $Y$ at $q$ is defined by:
\[
\cH[i]{q}(Y)=\Ext[\C^{\sf op}]{i}{S\{ q\}}{Y}
\quad \text{and} \quad
\hH[i]{q}(Y)=\Tor[\C]{i}{Y}{S\langle q\rangle}.
\]
\end{bfhpg}

\begin{lem}\label{homology}
Assume that $\C$ satisfies the strong retraction property. Let $X$ be an object in $\Mod{\C,A}$ and $I$ an injective cogenerator in $\Mod{\Bbbk}$. Then for each $i\geqslant 0$, there are isomorphisms
\[\cH[i]{q}(\Hom[\Bbbk]{X}{I})\is\Hom[\Bbbk]{\hH[i]{q}(X)}{I}
\quad
\text{and}
\quad
K'_q(\Hom[\Bbbk]{X}{I})\is\Hom[\Bbbk]{C_q(X)}{I}.\]
\end{lem}

\begin{prf*}
By definition, for each $i\geq 0$, we have
\begin{align*}
\cH[i]{q}(\Hom[\Bbbk]{X}{I})&=\Ext[\Cop]{i}{S\{ q\}}{\Hom[\Bbbk]{X}{I}}\\
&\cong\Hom[\Bbbk]{\Tor[\C]{i}{S\{q\}}{X}}{I}\\
&=\Hom[\Bbbk]{\hH[i]{q}(X)}{I},
\end{align*}
where the isomorphism holds by (\ref{1.9.1}) as $I$ is injective.

For the second isomorphism, we similarly have
\begin{align*}
K'_q(\Hom[\Bbbk]{X}{I})&=\Hom[\Cop]{S\{q\}}{\Hom[\Bbbk]{X}{I}}\\
&\cong\Hom[\Bbbk]{S\{q\}\otimes_\C X}{I}\\
&=\Hom[\Bbbk]{C_q(X)}{I},
\end{align*}
where the isomorphism follows from (\ref{1.9.1}).
\end{prf*}

The following result characterizes the flat objects in $\Mod{\C,A}$, complementing \cite[Theorem 7.29]{HJ21}.

\begin{thm}\label{flat2}
Assume that $\C$ is Hom-finite and satisfies the strong retraction property such that the pseudo-radical $\mathfrak{r}$ is nilpotent, that is, $\mathfrak{r}^N=0$ for some $N\in\mathbb{N}$. Let $X$ be an object in $\Mod{\C,A}$. Then $X$ is flat if and only if $\hH[1]{q}(X)=0$ and $C_q(X)$ is a flat $A$-module for each $q\in\C$.
\end{thm}
\begin{prf*}
Fix an injective cogenerator $I$ in $\Mod{\Bbbk}$. For the ``only if" part, assume that $X$ is flat. Then by Lemma \ref{flat-inj} one gets that $\Hom[\Bbbk]{X}{I}$ is injective in $\rMod{\C,A}$, and so it follows from \cite[Theorem 7.29(b)] {HJ21} that $\cH[1]{q}(\Hom[\Bbbk]{X}{I})=0$ and $K'_q(\Hom[\Bbbk]{X}{I})$ is an injective $\Aop$-module for each $q\in\C$.
By Lemma \ref{homology},
$$\cH[1]{q}(\Hom[\Bbbk]{X}{I})\is \Hom[\Bbbk]{\hH[1]{q}(X)}{I},$$
hence, $\hH[1]{q}(X)=0$. On the other hand, it follows from Lemma \ref{homology} that $$K'_q(\Hom[\Bbbk]{X}{I})\is\Hom[\Bbbk]{C_q(X)}{I},$$
so $C_q(X)$
is a flat $A$-module.

Conversely, assume that $\hH[1]{q}(X)=0$ and $C_q(X)$ is a flat $A$-module for each $q\in\C$. By Lemma \ref{homology}, we have $\cH[1]{q}(\Hom[\Bbbk]{X}{I})=0$ and $K'_q(\Hom[\Bbbk]{X}{I})$ is an injective $\Aop$-module. Hence, by \cite[Theorem 7.29(b)] {HJ21}, $\Hom[\Bbbk]{X}{I}$ is injective in $\rMod{\C,A}$. Applying Lemma \ref{flat-inj}, we conclude that $X$ is flat in $\Mod{\C,A}$.
\end{prf*}

\section{Construction of left Frobenius pairs}
\noindent
In this section, we develop a method for constructing left Frobenius pairs (equivalently, weak AB-contexts) in the category of $\C$-shaped diagrams, a useful tool for describing Gorenstein homological objects.

\begin{dfn}
For a subcategory $\X$ of $\Mod{A}$, define
\[\Fun(\C,\X)=\setof{X\in{\lMod{\C,A}}}{X(q)\in{\X}\ \mathrm{for\ each}\ q\in\C}\]
and
\[{\widetilde{\X}}=\Add\setof{G_q(N_q)}{q\in\C, N_q\in\X}.\]
\end{dfn}

\begin{lem}\label{preserve}
Assume that $\C$ is Hom-finite. Let $\X$ be a subcategory of $\Mod{A}$ closed under direct summands. Then for each $q\in\C$, we have ${G_q(\X)}\subseteq\Fun(\C,\X)$.
\end{lem}

\begin{prf*}
Fix $X\in \X$. For any $p,q\in \C$, we have $G_{q}(X)(p)=\Hom[\kk]{\C(p,q)}{X}$.
Since $\C$ is Hom-finite, $\C(p,q)$ is a finitely generated projective $\kk$-module.
But $\X$ is closed under finite coproducts and direct summands, it follows that $G_{q}(X)\in\Fun(\C,\X)$.
\end{prf*}

Recall from \cite[Definition 3.5]{HJ23} that a family $\{X_i\}_{i\in I}$ of objects in $\lMod{\C,A}$ is called \emph{locally finite} if for each $q\in\C$, the set $\setof{i\in I}{X_i(q)\neq 0}$ is finite. The next result follows from the proof of \cite[Proposition 3.7]{HJ23}.

\begin{lem}\label{locally finite}
Assume that $\C$ is locally bounded. Then for every family $\{M_q\}_{q\in\C}$ of objects in $\lMod{A}$, the family $\{G_q(M_q)\}_{q\in\C}$ is locally finite.
\end{lem}

\begin{lem}\label{coproduct}
Assume that $\C$ is Hom-finite and locally bounded. Let $\X$ be a subcategory of $\Mod{A}$ closed under direct summands. If $\{M_q\}_{q\in\C}$ is a family of objects in $\lMod{A}$ with each $M_q\in\X$, then $\prod_{q\in\C} G_q(M_q)\is\oplus_{q\in\C} G_q(M_q)$ lies in $\Fun(\C,\X)$.
\end{lem}

\begin{prf*}
The isomorphism holds by \cite[Proposition 3.7]{HJ23}.
For each $p\in\C$, the coproduct $\oplus_{q\in\C} (G_q(M_q))(p)$ is finite because the family $\{G_q(M_q)\}_{q\in\C}$ is locally finite by Lemma \ref{locally finite}.
Hence,
\[(\oplus_{q\in\C} G_q(M_q))(p) = \oplus_{q\in\C} (G_q(M_q))(p)\in\X\]
by Lemma \ref{preserve}, and therefore, $\oplus_{q\in\C} G_q(M_q) \in \Fun(\C,\X)$.
\end{prf*}

\begin{cor}\label{tilder}
Assume that $\C$ is Hom-finite and locally bounded. Let $\X$ be a subcategory of $\Mod{A}$ closed under direct summands. Then $\widetilde{\X}\subseteq\Fun(\C,\X)$.
\end{cor}

\begin{prf*}
Let $X\in\widetilde{\X}$. Then $X$ is a direct summand of an object $\oplus_{q\in\C}G_q(M_q)$ for some $M_q\in\X$. By Lemma \ref{coproduct}, this direct sum lies in $\Fun(\C,\X)$. Hence, $X(p)\in\X$ for each $p\in\C$.
\end{prf*}

\begin{lem}\label{cogenerator}
Assume that $\C$ is Hom-finite and locally bounded. Let $\X$ be a subcategory of $\Mod{A}$ closed under extensions and direct summands. If $\W$ is an injective cogenerator for $\X$, then $\widetilde{\W}$ is an injective cogenerator for $\Fun(\C,\X)$.
\end{lem}

\begin{prf*}
Fix an object $X\in \Fun(\C,\X)$. Then for each $q\in\C$ one has $X(q)\in\X$, and so there exists a short exact sequence
\[
0\to X(q)\to W_q \to T_q \to 0
\]
with $W_q\in\W$ and $T_q\in\X$. Since $G_q$ is exact, applying it gives a short exact sequence
\[
0\to \mathbb{G}(X)=\prod_{q\in\C} G_q(X(q))\to \prod_{q\in\C} G_q(W_q) \to \prod_{q\in\C} G_q(T_q) \to 0
\]
in $\lMod{\C,A}$, where $\prod_{q\in\C} G_q(W_q)\is\oplus_{q\in\C} G_q(W_q)\in\widetilde{\W}$ by \cite[Proposition 3.7]{HJ23} as $\C$ is locally bounded. On the other hand, it follows from \cite[Proposition 5.7]{HJ23} that there is an objectwise split short exact sequence
\[0 \to X \to \mathbb{G}(X) \to \mathbb{C}(X) \to 0.\]
Consider the push-out diagram
$$\xymatrix@C=15pt@R=15pt{&0\ar[d] &0\ar[d]\\
  & X \ar@{=}[r]\ar[d] & X \ar[d]\\
  0 \ar[r] & \mathbb{G}(X) \ar[r]\ar[d] & \prod_{q\in\C}G_q(W_q) \ar[r]\ar[d] & \prod_{q\in\C} G_q(T_q) \ar[r]\ar@{=}[d] & 0\\
  0 \ar[r] & \mathbb{C}(X) \ar[r]\ar[d] & Y \ar[r]\ar[d] & \prod_{q\in\C} G_q(T_q) \ar[r] & 0.\\
  & 0 & 0 }$$
Since $T_q, X(q) \in \X$, Lemma \ref{preserve} gives
$G_q(T_q), G_q(X(q))\in\Fun(\C,\X)$. By Lemma \ref{coproduct}, both $\prod_{q\in\C} G_q(T_q)\is\oplus_{q\in\C} G_q(T_q)$ and $\mathbb{G}(X)\is\oplus_{q\in\C} G_q(X(q))$ lie in $\Fun(\C,\X)$.
The first non-zero column in the above diagram is objectwise split exact, so, as $\X$ is closed under direct summands, $\mathbb{C}(X)\in\Fun(\C,\X)$.
Since $\Fun(\C,\X)$ is closed under extensions as so is $\X$ by the assumption, it follows that $Y\in\Fun(\C,\X)$. Consider the second non-zero column
\[0 \to X \to \prod_{q\in\C}G_q(W_q) \to Y \to 0\]
in the above diagram. Since $\prod_{q\in\C}G_q(W_q)\is \oplus_{q\in\C}G_q(W_q)\in \widetilde{\W}$ as $W_q\in\W$, one gets that $\widetilde{\W}$ is a cogenerator for $\Fun(\C,\X)$.

We then prove that $\Ext[\C,A]{i}{X}{W}=0$ for all $X\in\Fun(\C,\X)$, $W\in\widetilde{\W}$ and $i\geqslant1$. Without loss of generality, one may let
$$W=\oplus_{q\in\C} G_q(N_q)\is\prod_{q\in\C} G_q(N_q)$$
for some $N_q\in\W$. Thus it is enough to prove that $\Ext[\C,A]{i}{X}{G_q(N_q)}=0$. We mention that $(E_q, G_q)$ is an adjoint pair with $E_q$ and $G_q$ exact for each $q\in\C$. Then one has $\Ext[\C,A]{i}{X}{G_q(N_q)}\is\Ext[\C,A]{i}{X(q)}{N_q}=0$ by the assumption as $X(q)\in\X$. Thus, $\widetilde{\W}$ is an injective cogenerator for $\Fun(\C,\X)$.
\end{prf*}

\begin{thm}\label{Fpair}
Assume that $\C$ is Hom-finite and locally bounded. If $(\X, \W)$ is a left Frobenius pair in $\Mod{A}$, then $(\Fun(\C,\X), \widetilde{\W})$ is a left Frobenius pair in $\Mod{\C,A}$.
\end{thm}

\begin{prf*}
Since $(\X, \W)$ is a left Frobenius pair in $\Mod{A}$, one gets that $\X$ is closed under extensions, kernels of epimorphisms
and direct summands in $\Mod{A}$, and $\W$ is an injective cogenerator for $\X$ and closed under direct summands in $\Mod{A}$. Therefore, $\Fun(\C,\X)$ is closed under extensions, kernels of epimorphisms and direct summands in $\Mod{\C,A}$, and $\widetilde{\W}$ is an injective cogenerator for $\Fun(\C,\X)$ by Lemma \ref{cogenerator}. Moreover, ${\widetilde{\W}}$ is clearly closed under direct summands. Thus, $(\Fun(\C,\X), \widetilde{\W})$ is a left Frobenius pair in $\Mod{\C,A}$.
\end{prf*}

We now use Theorem \ref{Fpair} to construct a left Frobenius pair in $\Mod{\C,A}$ that will be needed later.

\begin{dfn}
Following \cite{SS20}, an object $X$ in $\Mod{\C, A}$ is called \emph{projectively coresolved Gorenstein flat object}, or a \emph{PGF object} for short, if there exists an exact sequence
\[\cdots \to P_{1}\to P_{0}\to P_{-1}\to\cdots\]
in $\Mod{\C, A}$ with each $P_{i}$ projective, such that $X\cong\Coker{(P_{1}\to P_{0})}$ and the sequence remains exact after applying $E\otimes_{\C,A}-$ for every injective object $E$ in $\rMod{\C,A}$. Let $\PGF{\C,A}$ denote the subcategory of PGF objects in $\Mod{\C,A}$. When $A=\kk$, we write $\PGF{\C}$ instead of $\PGF{\C,A}$.
\end{dfn}

Let $\PGF{A}$ denote the subcategory of PGF $A$-modules, and $\Prj{A}$ the subcategory of projective $A$-modules.

\begin{cor}\label{Fpair-exa}
Assume that $\C$ is Hom-finite and locally bounded. Then the pair $(\Fun(\C,\PGF{A}), \widetilde{\Prj{A}})$ is a left Frobenius pair in $\Mod{\C,A}$.
\end{cor}

\begin{prf*}
Since $(\PGF{A}, \Prj{A})$ is a left Frobenius pair in $\Mod{A}$ by \cite[Example 2.25]{LY22}, Theorem \ref{Fpair} yields the desired conclusion.
\end{prf*}

Let $\Flat{A}$ denote the subcategory of flat $A$-modules, and $\FC{A}$ the subcategory of flat and cotorsion $A$-modules.

\begin{cor}\label{Fpair-flat}
Assume that $\C$ is Hom-finite and locally bounded. Then the pair $(\Fun(\C,\Flat{A}), \widetilde{\FC{A}})$ is a left Frobenius pair in $\Mod{\C,A}$.
\end{cor}

\begin{prf*}
It is known that $(\Flat{A}, \FC{A})$ is a left Frobenius pair in $\Mod{A}$ and $\Flat{A}$ is closed under coproducts and direct summands; see \cite[Example 2.20]{LY22}. Hence, by Theorem \ref{Fpair}, $(\Fun(\C,\Flat{A}), \widetilde{\FC{A}})$ is a left Frobenius pair in $\Mod{\C,A}$.
\end{prf*}

We close this section with a partial converse to Theorem \ref{Fpair}.

\begin{prp}
Assume that $\C$ is Hom-finite and locally bounded, and satisfies the strong retraction property. Let $\X$ and $\W$ be subcategories of $\Mod{A}$ with $\W$ closed under direct summands. Then the following statements are equivalent.
\begin{eqc}
\item $(\X, \W)$ is a left Frobenius pair in $\Mod{A}$.
\item $(\Fun(\C,\X), \widetilde{\W})$ is a left Frobenius pair in $\Mod{\C,A}$.
\end{eqc}
\end{prp}

\begin{prf*}
\eqclbl{i}$\implies$\eqclbl{ii} follows from Proposition \ref{Fpair}.

\proofofimp{ii}{i} Let $0 \to K \to M \to C \to 0$ be an exact sequence in $\Mod{A}$ with $C\in\X$. Applying $S_q$ gives an exact sequence
$0 \to S_q(K) \to S_q(M) \to S_q(C) \to 0$ in $\Mod{\C,A}$.
Since $S_q(C)$ is in $\Fun(\C,\X)$ and $\Fun(\C,\X)$ is closed under extensions and kernels of epimorphisms, it follows that
\[
K \in \X \, \Leftrightarrow \, S_q(K) \in \Fun(\C,\X) \, \Leftrightarrow \, S_q(M) \in \Fun(\C,\X) \, \Leftrightarrow \, M \in \X.
\]
Thus $\X$ is closed under extensions and kernels of epimorphisms. Similarly, one may prove that $\X$ is closed under direct summands as so is $\Fun(\C,\X)$ via the functor $S_q$.

We next show that $\W$ is an injective cogenerator for $\X$.
Fix an object $N$ in $\X$. Then $S_q(N)\in \Fun(\C,\X)$ for each $q\in\C$, and so there is an exact sequence
$$0 \to S_q(N) \to W \to X \to 0$$
in $\Mod{\C,A}$ with $W\in\widetilde{\W}$ and $X\in\Fun(\C,\X)$. This yields an exact sequence
$$0 \to N \to W(q) \to X(q) \to 0$$
with $X(q)\in\X$ and $W(q)\in\W$ by Corollary \ref{tilder}. Hence, $\W$ is a cogenerator for $\X$.

Finally, for each $Y\in\X$, $V\in\W$ and $i\geqslant1$, one has
\begin{align*}
\Ext[A]{i}{Y}{V}&=\Ext[A]{i}{E_qS_q(Y)}{V}\\
&\cong\Ext[\C,A]{i}{S_q(Y)}{G_q(V)}\\
&=0,
\end{align*}
where the first isomorphism holds as $(E_q, G_q)$ is an adjoint pair with both $E_q$ and $G_q$ exact for each $q\in\C$, and the last equality holds as $S_q(Y)\in\Fun(\C,\X)$ and $G_q(V)\in\widetilde{\W}$. Thus, $\W$ is an injective cogenerator for $\X$, so $(\X,\W)$ is a left Frobenius pair in $\Mod{A}$.
\end{prf*}

\section{A description for PGF objects}
\noindent
In this section, we provide a description of PGF objects in $\lMod{\C,A}$ using the results established in the previous section.

\begin{lem}\label{containment}
Assume that $\C$ is Hom-finite, locally bounded and has a Serre functor $\mathbb{S}$. If $X$ is a PGF object in $\lMod{\C,A}$, then $X(q)$ is a PGF $A$-module for each $q\in\C$.
\end{lem}

\begin{prf*}
By definition, there is an exact sequence
\[P=\cdots \to P_{1}\to P_{0}\to P_{-1}\to\cdots\]
of projective objects in $\lMod{\C,A}$ such that $X\cong\Coker{(P_{1}\to P_{0})}$ and the sequence remains exact after applying $E\otimes_{\C,A}-$ for every injective object $E$ in $\rMod{\C,A}$. Set $X_i=\Coker{(P_{i+1}\to P_{i})}$. Then  $\Tor[1]{\C,A}{E}{X_i}=0$. For each $q\in\C$, one gets an exact sequence
\[P(q)=\cdots \to P_{1}(q)\to P_{0}(q)\to P_{-1}(q)\to\cdots\]
of $A$-modules with each $P_{i}(q)$ projective, since the functor $E_q$ preserves projectives. By Lemma \ref{Tor vanish}, $\Tor[1]{A}{I}{X_i(q)}=0$ for each injective module $I\in\rMod{A}$. Hence, the sequence $P(q)$ remains exact after applying $I\otimes_{A}-$, so $X(q)$ is a PGF $A$-module for each $q\in\C$.
\end{prf*}

\begin{thm}\label{objectwise}
Assume that $\C$ is Hom-finite, locally bounded and has a Serre functor $\mathbb{S}$. Then $\PGF{\C,A}=\Fun(\C,\PGF{A})$.
\end{thm}

\begin{prf*}
The inclusion $\PGF{\C,A}\subseteq\Fun(\C,\PGF{A})$ holds by Lemma \ref{containment}. It remains to prove the reverse inclusion. Take $X\in\Fun(\C,\PGF{A})$. By Corollary \ref{Fpair-exa}, there is an exact sequence
\begin{equation*}\label{diag01}
\tag{\ref{objectwise}.1}
0 \to X \to P_{-1} \to X_{-1} \to 0
\end{equation*}
with $P_{-1}\in\widetilde{\Prj{A}}$ and $X_{-1}\in\Fun(\C,\PGF{A})$.
By \cite[Lemma 3.4]{HJ23},
\[\widetilde{\Prj{A}}\is\Add\setof{F_q(N_q)}{q\in\C, N_q\in\Prj{A}},\]
and by \cite[Proposition 3.12]{HJ21}, $\widetilde{\Prj{A}}$ coincides with the subcategory of projective objects in $\lMod{\C,A}$.
Thus, for each $q\in\C$, the sequence
\[0 \to X(q) \to P_{-1}(q) \to X_{-1}(q) \to 0\]
is exact, where $P_{-1}(q)$ is projective and $X_{-1}(q)$ is PGF by Lemma \ref{containment}. So one has $\Tor[1]{A}{I}{X_{-1}(q)}=0$ for each injective $\Aop$-module $I$, which yields that
\[\Tor[1]{\C,A}{E}{X_{-1}}=0\]
for each injective object $E$ in $\rMod{\C,A}$ by Lemma \ref{Tor vanish}.
Therefore, the sequence (\ref{diag01}) remains exact after applying $E\otimes_{\C,A}-$. Repeating this process yields an exact sequence
\begin{equation*}\label{diag02}
\tag{\ref{objectwise}.2}
0 \to X \to P_{-1} \to P_{-2} \to \cdots
\end{equation*}
in $\lMod{\C,A}$ with each $P_i$ projective such that it remains exact after applying $E\otimes_{\C,A}-$ for each injective object $E$ in $\rMod{\C,A}$.

Now take a projective resolution
\begin{equation*}\label{diag03}
\tag{\ref{objectwise}.3}
\cdots \to P_1 \to P_0 \to X \to 0.
\end{equation*}
Since $X(q)$ is a PGF $A$-module for each $q\in\C$, one has $\Tor[i]{A}{I}{X(q)}=0$ for each injective $\Aop$-module $I$ and each $i\geq 1$. By Lemma \ref{Tor vanish}, it follows that $\Tor[i]{\C,A}{E}{X}=0$ for each injective object $E\in\rMod{\C,A}$. Hence, the sequence (\ref{diag03}) remains exact after applying $E\otimes_{\C,A}-$. Combining (\ref{diag02}) and (\ref{diag03}), one concludes that $X$ is a PGF object in $\lMod{\C,A}$, that is, $X\in\PGF{\C,A}$.
\end{prf*}

\section{A flat model structure on $\Mod{\C,A}$}
\noindent
In this section, we construct a flat model structure on $\lMod{\C,A}$ using PGF objects. Let $G$ be a PGF object in $\Mod{\C,A}$.
By definition, there is an exact sequence
\[\cdots \to P_{1}\xra{d_{1}} P_{0}\xra{d_{0}} P_{-1}\to\cdots\]
in $\Mod{\C, A}$ with each $P_{i}$ projective such that $X\cong\Coker{d_1}$ and the sequence remains exact after applying $E\otimes_{\C,A}-$ for every injective object $E$ in $\rMod{\C,A}$.
\textbf{In this section}, denote by $\Omega^{i}G=\mathrm{Im} d_{i}$.
It is routine to verify that $\Omega^{i}G$ is unique up to projective summands.
Hence, the following definition is well-defined.

\begin{dfn}
Let $\G$ be a class of PGF objects in $\Mod{\C,A}$. Define
\[\W_{\G} = \{X\in \Mod{\C,A}\mid \Ext[\C,A]{1}{\Omega^{n}G}{X}=0 \ \text{for all} \ G\in \G \ \text{and}\ n\in \mathbb{Z}\}.\]
\end{dfn}

\begin{rmk}\label{generate}
It is clear that $(^{\bot}\WG,\WG)$ is a cotorsion pair generated by the class $\{\Omega^{n}G\mid G\in \G, n\in \mathbb{Z}\}.$
\end{rmk}

Recall from \cite{rha} that an object $X$ in $\Mod{\C, A}$ is called \emph{Gorenstein projective} if there is an exact sequence
\[\cdots \to P_{1}\to P_{0}\to P_{-1} \to \cdots\]
in $\Mod{\C,A}$ with each $P_{i}$ projective such that $X\cong\Coker{(P_{1}\to P_{0})}$ and the sequence remains exact after applying $\Hom[\C,A]{-}{P}$ for every projective object $P$ in $\Mod{\C,A}$. The following result is known for modules over a ring; see \cite[Theorem 4.4]{SS20}.

\begin{lem}\label{lem:PGF-GP}
Assume that $\C$ is Hom-finite, locally bounded and has a Serre functor $\mathbb{S}$. Let $X$ be a PGF object in $\Mod{\C,A}$. Then $X$ is Gorenstein projective in $\Mod{\C,A}.$
\end{lem}

\begin{prf*}
Since $X$ is PGF in $\Mod{\C,A}$, Theorem \ref{objectwise} implies that $X(q)$ is a PGF $A$-module for each $q\in\C$. By \cite[Theorem 4.4]{SS20}, each $X(q)$ is a Gorenstein projective $A$-module.
It then follows from \cite[Theorem B]{DLLM} that $X$ is Gorenstein projective in $\Mod{\C,A}$.
\end{prf*}

\begin{prp}\label{prop:PGF-cotorsion pair}
Assume that $\C$ is Hom-finite, locally bounded, and has a Serre functor $\mathbb{S}$. Let $\G$ be a class of $PGF$ objects in $\Mod{\C,A}$.
Then the following statements hold:
\begin{prt}
\item The subcategory $\WG$ is thick, and $^{\bot}\WG\cap\WG=\Prj{\C,A}.$
\item If $\G$ is a set, then $(^{\bot}\WG,\WG)$ is a hereditary complete cotorsion
      pair.
\end{prt}
\end{prp}

\begin{prf*}
(a). By definition, $\WG$ is closed under direct summands.
Consider exact sequence
\[0\to X_{1}\to X_{2}\to X_{3}\to 0\]
in $\Mod{\C,A}$.
For any $G\in \G$ and $n\in \mathbb{Z}$, we have the following exact sequence
\[\Ext[\C,A]{1}{\Omega^{n}G}{X_{1}}\to \Ext[\C,A]{1}{\Omega^{n}G}{X_{2}}\to \Ext[\C,A]{1}{\Omega^{n}G}{X_{3}}\to \Ext[\C,A]{2}{\Omega^{n}G}{X_{1}}.\]
We mention that $\Ext[\C,A]{2}{\Omega^{n}G}{X_{1}}\cong \Ext[\C,A]{1}{\Omega^{n+1}G}{X_{1}}$.
It follows that whenever $X_{1}$ and one of $X_{2}$ or $X_{3}$ lie in $\WG$, then so does the third of them.
Now assume that $X_{2}, X_{3}\in \WG$.
The above sequence gives $\Ext[\C,A]{2}{\Omega^{n}G}{X_{1}}=0$ for all $n$, and hence
\[\Ext[\C,A]{1}{\Omega^{n}G}{X_{1}}\cong\Ext[\C,A]{2}{\Omega^{n-1}G}{X_{1}}=0,\]
so $X_1\in\WG$. Thus, $\WG$ is thick.

Next, we prove $\Prj{\C,A}=\WG\cap {^{\bot}\WG}$. Clearly, $\Prj{\C,A}\subseteq{^{\bot}\WG}$, and by Lemma \ref{lem:PGF-GP}, $\Prj{\C,A}\subseteq\WG$, so $\Prj{\C,A}\subseteq\WG\cap {^{\bot}\WG}$. Conversely, take $X\in \WG\cap {^{\bot}\WG}$ and consider an exact sequence
$$0\to K\to P\to X\to 0$$
with $P\in\Prj{\C,A}$. As above, $P\in\WG$, and since $\WG$ is thick and $X\in\WG$, we get $K\in\WG$. But $X\in{^{\bot}\WG}$, so $\Ext[\C,A]{1}{X}{K}=0$. Hence, the sequence splits, giving $X\in\Prj{\C,A}$.

(b). By Remark \ref{generate}, $(^{\bot}\WG,\WG)$ is a cotorsion pair generated by
\[\{\Omega^{n}G\mid G\in\G,\ n\in\mathbb{Z}\},\]
which is a set whenever $\G$ is. Thus, $(^{\bot}\WG,\WG)$ is complete by, for example, \cite[Theorem 6.11(b)]{RT12}. Moreover, since $\WG$ is thick, the cotorsion pair is hereditary.
\end{prf*}

\begin{lem}\label{flat cotorsion pair}
Assume that $\C$ is Hom-finite, locally bounded, and has a Serre functor $\mathbb{S}$. Let $\G$ be a set of $PGF$ objects in $\Mod{\C,A}$. Then the pair
\[(^{\perp}(\Cot{\C,A}\cap\WG), \, \Cot{\C,A}\cap\WG)\]
is a complete and hereditary cotorsion pair in $\Mod{\C,A}$ generated by a set, and $^{\perp}(\Cot{\C,A}\cap\WG)\cap\WG=\Flat{\C,A}$.
\end{lem}

\begin{prf*}
It is known that the cotorsion pair $(\Flat{\C,A}, \Cot{\C,A})$ is generated by a set $\Se$ (see \ref{representable}), that is, $\Se \subseteq \Flat{\C,A}$ with $\Se^\perp = \Cot{\C,A}$. Moreover, by Remark \ref{generate},
$\WG=\{\Omega^{n}G\mid G\in\G,\ n\in\mathbb{Z}\}^{\perp}$.
Since $\G$ is a set, the generating class
\[\{\Omega^{n}G\mid G\in\G,\ n\in\mathbb{Z}\}\]
is also a set.
Hence,
\[(\Se \cup \{\Omega^{n}G\mid G\in\G,n\in\mathbb{Z}\})^\perp = \Se^\perp \cap \{\Omega^{n}G\mid G\in\G,n\in\mathbb{Z}\}^\perp = \Cot{\C,A} \cap \WG,\]
so $(^{\perp}(\Cot{\C,A} \cap \WG), \Cot{\C,A} \cap \WG)$ is a complete cotorsion pair in $\Mod{\C,A}$ generated by a set; see \cite[Theorem 6.11(b)]{RT12}.
Since both $\Cot{\C,A}$ and $\WG$ are closed under cokernels of monomorphisms (the latter by Proposition \ref{prop:PGF-cotorsion pair}), so is $\Cot{\C,A} \cap \WG$. Hence, the cotorsion pair is hereditary.

It remains to prove the equality
\[
^{\perp}(\Cot{\C,A} \cap \WG) \cap \WG = \Flat{\C,A}.
\]
Clearly,
$\Flat{\C,A} = {^{\perp}\Cot{\C,A}} \subseteq {^{\perp}(\Cot{\C,A} \cap \WG)}$.
Moreover, if $X \in \Flat{\C,A}$, then $X$ is in $\WG$ by \cite[Theorem 4.4]{SS20}. Thus we have
$$\Flat{\C,A} \subseteq{^{\perp} (\Cot{\C,A} \cap \WG) \cap \WG}.$$
For the reverse inclusion, take $Y \in{^{\perp}(\Cot{\C,A} \cap \WG) \cap\WG}$. Since the cotorsion pair $(\Flat{\C,A}, \Cot{\C,A})$ is complete (see \ref{representable}), there is an exact sequence
\begin{equation*}\label{eq:ipg1}\tag{\ref{flat cotorsion pair}.1}
0 \to C \to F \to Y \to 0
\end{equation*}
with $C \in \Cot{\C,A}$ and $F \in \Flat{\C,A} \subseteq \WG$. Since $\WG$ is thick, $C\in\WG$. Hence, $\Ext[\C,A]{1}{Y}{C}=0$ because $Y\in{^{\perp}(\Cot{\C,A}\cap\WG)}$. Thus, the sequence splits, giving $Y\in\Flat{\C,A}$.
\end{prf*}

We are now ready to construct a hereditary abelian model structure in $\Mod{\C,A}$.

\begin{thm}\label{thm:model structure-PGF}
Assume that $\C$ is Hom-finite, locally bounded, and has a Serre functor $\mathbb{S}$. Let $\G$ be a set of PGF objects in $\Mod{\C,A}$.
Then there exists a hereditary abelian model structure
\[({^{\perp}(\Cot{\C,A}\cap\WG)}, \, \WG, \, \Cot{\C,A})\]
in $\Mod{\C,A}$.
\end{thm}

\begin{prf*}
It follows from Lemma \ref{flat cotorsion pair} that
\[
(^{\perp}(\Cot{\C,A}\cap\WG), \, \Cot{\C,A}\cap\WG)
\]
is a complete and hereditary cotorsion pair in $\Mod{\C,A}$, and
\[
^{\perp}(\Cot{\C,A}\cap\WG)\cap\WG=\Flat{\C,A}.
\]
Thus, $(^{\perp}(\Cot{\C,A}\cap\WG)\cap\WG, \Cot{\C,A})$ is a complete and hereditary cotorsion pair in $\Mod{\C,A}$; see \ref{representable}. By Proposition \ref{prop:PGF-cotorsion pair}, $\WG$ is a thick subcategory. Hence, $({^{\perp}(\Cot{\C,A}\cap\WG)}, \WG, \Cot{\C,A})$ is a hereditary abelian model structure in $\Mod{\C,A}$.
\end{prf*}

\section{Applications to Q-shaped derived categories}
\noindent
In this section, we apply Theorem \ref{thm:model structure-PGF} to some special cases and give flat models for $\Q$-shaped derived categories, and describe the cofibrant objects. The results in this section require the following setup from \cite{HJ21}.

\begin{setup}
Throughout this section, assume that the ring $\kk$ is noetherian and hereditary (for example, $\kk$ is a field or $\kk=\mathbb{Z}$), and let $\Q$ be a small $\kk$-linear category satisfying:
\begin{rqm}
\item[$\bullet$]
     $\Q$ is Hom-finite, locally bounded, has a Serre functor $\mathbb{S}$, and satisfies the strong retraction property.
\item[$\bullet$]
     The pseudo-radical $\mathfrak{r}$ is nilpotent, that is, $\mathfrak{r}^N=0$ for some $N\in\mathbb{N}$.
\end{rqm}
\end{setup}

Fix $\G=\{S_{q}(A)\mid q\in \Q\}$. By Theorem \ref{objectwise}, $\G$ is a set of PGF objects in $\lMod{\Q,A}$. We now describe the trivial objects in the abelian model structure from Theorem \ref{thm:model structure-PGF}.

\begin{lem}\label{wg}
Let $\G=\{S_{q}(A)\mid q\in \Q\}$ and $X$ an object in $\lMod{\Q,A}$. Then $X$ lies in $\WG$ if and only if $\cH[1]{q}(X)=0$ for each $q\in\Q$.
\end{lem}

\begin{prf*}
It follows from \cite[Lemma 5.1]{SV26} that
$$\WG=\{X\in \Mod{\Q,A}\mid \Ext[\C,A]{1}{S_{q}(A)}{X}=0 \ \text{for all} \ q\in \Q \}.$$
Indeed, we have the following isomorphisms:
\begin{align*}
\Ext[\Q,A]{1}{S_{q}(A)}{X}&=\Ext[\Q,A]{1}{S\langle q\rangle\otimes_{\kk}A}{X}\\
&\cong\Ext[\Q]{1}{S\langle q\rangle}{X^{\natural}}\\
&= \mathbb{H}_{[q]}^{1}(X^{\natural})\\
&= \mathbb{H}_{[q]}^{1}(X)^{\natural},
\end{align*}
where the isomorphism follows from Lemma \ref{iso}(a) as $S\langle q\rangle$ is a PGF object in $\Mod{\Q}$ by Theorem \ref{objectwise}, and
the last equality holds by \cite[Remark 7.12]{HJ21}.
Thus, $\Ext[\Q,A]{1}{S_{q}(A)}{X}$ vanishes exactly when $\mathbb{H}_{[q]}^{1}(X)$ does. Consequently, $X\in\WG$ if and only if $\cH[1]{q}(X)=0$ for each $q\in\Q$.
\end{prf*}

Following \cite[Definition 4.1]{HJ21}, let $\E$ denote the class of \emph{exact objects} in $\lMod{\Q,A}$, that is,
$$\E=\setof{X\in\lMod{\Q,A}}{\pd_{\Q}X^{\natural}<\infty\ \mathrm{or}\ \id_{\Q}X^{\natural}<\infty}.$$

The following result, which can be found in \cite[Theorem 7.1]{HJ21}, is used frequently in the rest of the paper.

\begin{lem}\label{thm7.1}
Let $X$ be an object in $\lMod{\Q,A}$. Then the following statements are equivalent.
\begin{eqc}
\item $X\in\E$.
\item $\cH[i]{q}{(X)}=0$ for each $q\in\Q$ and $i>0$.
\item $\cH[1]{q}{(X)}=0$ for each $q\in\Q$.
\item $\hH[i]{q}{(X)}=0$ for each $q\in\Q$ and $i>0$.
\item $\hH[1]{q}{(X)}=0$ for each $q\in\Q$.
\end{eqc}
\end{lem}

We obtain a flat model structure whose homotopy category is the $\Q$-shaped derived category introduced in \cite{HJ21}; this model structure was first given in \cite[Corollary 2.8]{DLLM}.

\begin{thm}\label{flat model structure}
There is a hereditary abelian model structure
$$(^{\perp}(\Cot{\Q,A}\cap\E), \, \E, \, \Cot{\Q,A})$$
on $\lMod{\Q,A}$, whose homotopy category is the $\Q$-shaped derived category.
\end{thm}

\begin{prf*}
Let $\G=\{S_{q}(A)\mid q\in \Q\}$. Then by Theorem \ref{objectwise}, $\G$ is a set of PGF objects in $\lMod{\Q,A}$, and $\WG=\E$ by Lemmas \ref{wg} and \ref{thm7.1}. Thus by Theorem \ref{thm:model structure-PGF},
\[(^{\perp}(\Cot{\Q,A}\cap\E), \E, \Cot{\Q,A})\]
is a hereditary abelian model structure on $\lMod{\Q,A}$, whose trivial objects coincide with those in the projective (injective) model structure given in \cite[Theorem A]{HJ21}. Thus the corresponding homotopy category is the $\Q$-shaped derived category.
\end{prf*}

We now characterize the cofibrant objects in $^{\perp}(\Cot{\Q,A}\cap\E)$.

\begin{lem}\label{dgflat}
Let $X$ be an object in $^\perp(\Cot{\Q,A}\cap\E)$ and $I$ an injective cogenerator in $\Mod{\Bbbk}$. Then $\Hom[\Bbbk]{X}{I}$ lies in $\E^\perp$.
\end{lem}

\begin{prf*}
Take $X\in {}^\perp(\Cot{\Q,A}\cap\E)$, and let $E\in\rMod{\Q,A}$ with $E\in\E$.
Then by Lemma \ref{thm7.1}, $\hH[1]{q}(E)=0$ for each $q\in\Q$.
It follows from (\ref{1.9.1}) that
\[
\Ext[\Q^{\sf op},A^{\sf op}]{1}{E}{\Hom[\Bbbk]{X}{I}}\is\Ext[\Q,A]{1}{X}{\Hom[\Bbbk]{E}{I}}.\]
By Lemma \ref{homology},
$$\cH[1]{q}(\Hom[\Bbbk]{E}{I})\is\Hom[\Bbbk]{\hH[1]{q}(E)}{I}=0$$
for each $q\in\Q$. Hence, $\Hom[\Bbbk]{E}{I}\in\E$ by Lemma \ref{thm7.1}.
Moreover, by Corollary \ref{cotorsion}, $\Hom[\Bbbk]{E}{I}\in\Cot{\Q,A}$, so $\Ext[\Q,A]{1}{X}{\Hom[\Bbbk]{E}{I}}=0$ because $X\in{^\perp(\Cot{\Q,A}\cap\E)}$. Therefore,
\[\Ext[\Q^{\sf op},A^{\sf op}]{1}{E}{\Hom[\Bbbk]{X}{I}}=0,\]
hence, $\Hom[\Bbbk]{X}{I}\in\E^\perp$.
\end{prf*}

The following result is taken from \cite[Theorem E]{HJ23}; we include a proof for the reader's convenience.

\begin{lem}\label{inj}
If $X\in\E^\perp$, then $X(q)$ is injective for each $q\in\Q$.
\end{lem}

\begin{prf*}
It follows from \ref{adjoint triple} that $(F_q, E_q)$ is an adjoint pair for each $q\in\Q$ with both $F_q$ and $E_q$ exact. Hence, for each $N\in\Mod{A}$, one has \[\Ext[A]{1}{N}{E_q(X)}\is\Ext[\Q,A]{1}{F_q(N)}{X}=0,\]
since $F_q(N)\in\E$ by \cite[Lemma 7.14]{HJ21} and Lemma \ref{thm7.1}.
\end{prf*}

\begin{prp}\label{cofibrant}
Let $X$ be an object in $^\perp(\Cot{\Q,A}\cap\E)$. Then each $A$-module $X(q)$ is flat for $q\in\Q$.
\end{prp}

\begin{prf*}
Let $I$ be an injective cogenerator in $\Mod{\Bbbk}$. Then by Lemma \ref{dgflat}, $\Hom[\Bbbk]{X}{I}\in\E^\perp$. Lemma \ref{inj} then gives that  $\Hom[\Bbbk]{X(q)}{I}=\Hom[\Bbbk]{X}{I}(q)$ is an injective $\Aop$-module for each $q\in\Q$, and so $X(q)$ is flat.
\end{prf*}

\begin{lem}\label{Add}
Let $\X$ be a subcategory of $\Mod{A}$ closed under direct summands, and $X$ in ${\widetilde{\X}}$. Then $X\in\E$ and $C_q(X)\in\X$ for each $q\in\Q$.
\end{lem}

\begin{prf*}
Since $X$ is in ${\widetilde{\X}}$, it is a direct summand of $\oplus_{p\in\C}G_p(N_p)$ for some $N_p\in\X$. By \cite[Proposition 3.7]{HJ23},
\[
\oplus_{p\in\C}G_p(N_p)\is\prod_{p\in\C}G_p(N_p).
\]
Then \cite[Lemma 7.14]{HJ21} gives $\cH[1]{q}(\oplus_{p\in\C}G_p(N_p))=0$ for each $q\in\Q$, so one has $\oplus_{p\in\C}G_p(N_p)\in\E$ by Lemma \ref{thm7.1}, and hence $X\in\E$.

On the other hand, for each $q\in\Q$, $C_q(X)$ is a direct summand of $$C_q(\oplus_{p\in\C}G_p(N_p))\is\oplus_{p\in\C}C_qG_p(N_p)\is \oplus_{p\in\C}C_qF_{\mathbb{S}^{-1}(p)}(N_p)=N_{\mathbb{S}(q)}\in\X,$$
where the first isomorphism holds as $(C_q,S_q)$ is an adjoint pair, the second one follows from \cite[Lemma 3.4]{HJ23}, and the equality holds by \cite[Lemma 7.28(a)]{HJ21}. Since $\X$ is closed under direct summands, one gets that $C_q(X)\in\X$ for each $q\in\Q$.
\end{prf*}

\begin{thm}\label{cofibrant object}
Assume that $A$ has finite weak global dimension. Let $X$ be an object in $\Mod{\C,A}$. Then $X\in{^\perp(\Cot{\Q,A}\cap\E)}$ if and only if $X(q)$ is flat in $\Mod{A}$ for every $q\in\Q$.
\end{thm}

\begin{prf*}
The ``only if" part holds by Proposition \ref{cofibrant}. For the ``if" part, we assume that the weak global dimension of $A$ is $n$, and that $X(q)$ is flat for each $q\in\Q$. It follows from Theorem \ref{flat model structure} that $$({^\perp(\Cot{\Q,A}\cap\E)},\Cot{\Q,A}\cap\E)$$
is a hereditary complete cotorsion pair. So there is a short exact sequence
\begin{equation*}\label{4.11.1}
\tag{\ref{cofibrant object}.1}
0\to K\to F \to X \to 0
\end{equation*}
in $\Mod{\Q,A}$ with $F\in{^\perp(\Cot{\Q,A}\cap\E)}$ and $K\in\Cot{\Q,A}\cap\E$, which yields a short exact sequence
$$0 \to K(q) \to F(q) \to X(q) \to 0$$
 of $A$-modules for each $q\in\Q$. Since $X(q)$ and $F(q)$ are flat by Proposition \ref{cofibrant}, one has $K(q)$ is flat, and so $K$ is in $\Fun(\Q,\Flat{A})$. We mention that $\widetilde{\FC{A}}$ is an injective cogenerator for $\Fun(\Q,\Flat{A})$ as $(\Fun(\Q,\Flat{A}), \widetilde{\FC{A}})$ is a left Frobenius pair in $\Mod{\Q,A}$ by Corollary \ref{Fpair-flat}. Then there is an exact sequence
$$0 \to K \to P^0 \to \cdots \to P^{n-1} \to K' \to 0$$
in $\Mod{\Q,A}$ with each $P^i\in\widetilde{\FC{A}}$ and $K'\in\Fun(\Q,\Flat{A})$. Note that $K$ is in $\E$ and each $P^i$ is in $\E$ by Lemma \ref{Add}. Then $K'$ is in $\E$ as $\E$ is a thick subcategory. Thus one has $L_iC_q(K')=\hH[i]{q}(K')=0$ for each $i\geqslant1$ by Lemma \ref{thm7.1}. Consequently, the sequence
$$0 \to C_q(K) \to C_q(P^0) \to \cdots \to C_q(P^{n-1}) \to C_q(K') \to 0$$
is exact, where each $C_q(P^i)$ is flat by Lemma \ref{Add}. Since $A$ has finite weak global dimension $n$, one gets that $C_q(K)$ is flat. On the other hand, one has $\hH[1]{q}(K)=0$ by Lemma \ref{thm7.1} as $K\in\E$, and so $K$ is flat by Theorem \ref{flat2}. Thus, $K$ is in $\Flat{\C,A}\cap\Cot{\C,A}$.

Since all flat $A$-modules are Gorenstein flat, $X$ is Gorenstein flat by \cite[Theorem 3.11]{DLLM}. Then \cite[Theorem 4.11]{SS20} yields $\Ext[\Q,A]{1}{X}{K}=0$, so the sequence (\ref{4.11.1}) splits. Therefore, $X$ is a direct summand of $F\in{^\perp(\Cot{\Q,A}\cap\E)}$, and hence, $X\in{^\perp(\Cot{\Q,A}\cap\E)}$.
\end{prf*}

\begin{lem}\label{lem4.13}
The following statements are equivalent.
\begin{eqc}
\item For every $C\in\Cot{\Q,A}\cap\E$ and every $q\in \Q$, $K_q(C)$ is cotorsion.
\item For every flat $A$-module $F$ and every $q\in \Q$, $S_q(F)$ lies in
      $^\perp(\Cot{\Q,A}\cap\E)$.
\end{eqc}
\end{lem}

\begin{prf*}
\proofofimp{i}{ii} Let $C$ be in $\Cot{\Q,A}\cap\E$. Then $K_q(C)$ is cotorsion for all $q\in\Q$. Since $(S_q,K_q)$ is an adjoint pair with $S_q$ exact for each $q\in\C$, one has
\[
\Ext[\C,A]{1}{S_q(F)}{C} \is \Ext[A]{1}{F}{K_q(C)} = 0
\]
for each flat $A$-module $F$, where the isomorphism holds by \cite[Lemma 1.3]{HJ19b} as $R^1K_q(C)=\Ext[\C]{1}{S\langle q\rangle}{C}=\cH[1]{q}(C)=0$. Thus, one has $S_q(F)\in{^\perp(\Cot{\C,A}\cap\E)}$.

\proofofimp{ii}{i} Let $F$ be a flat $A$-module. Then $S_q(F)$ is in $^\perp(\Cot{\Q,A}\cap\E)$ for each $q\in\Q$, so
\[
\Ext[A]{1}{F}{K_q(C)} \is \Ext[\C,A]{1}{S_q(F)}{C} = 0
\]
for each $C\in\Cot{\Q,A}\cap\E$, where the isomorphism holds as proved before. This implies that $K_q(C)$ is cotorsion for each $q\in\Q$.
\end{prf*}

The following is an immediate consequence of Theorem \ref{cofibrant object} and Lemma \ref{lem4.13}.

\begin{cor}\label{dg-cotorsion}
Assume that $A$ has finite weak global dimension. Then for every $C\in\Cot{\Q,A}\cap\E$ and every $q\in\Q$, the module $K_q(C)$ is cotorsion.
\end{cor}

\begin{rmk}
Consider the quiver
$$\Gamma=\cdots \to \bullet \xra{\partial} \bullet \xra{\partial} \bullet \to \cdots$$
with the relations $\partial^2=0$, and let $\Q$ be its path category. Then $\lMod{\Q,A}$ identifies with the category of chain complexes $\mathsf{Ch}(A)$.
In this setting, $\Cot{\Q,A}$ is the subcategory of dg-cotorsion complexes, and $\E$ is the subcategory of exact complexes. It is known that every exact dg-cotorsion complex has cotorsion cycles (\cite[Theorem 1.3]{BCE20}). However, it remains unclear whether this property extends to the more general context of $\Q$-shaped diagrams. More precisely, we do not know whether Corollary \ref{dg-cotorsion} holds without the finite weak global dimension assumption on $A$.
\end{rmk}


\bibliographystyle{amsplain-nodash}

\begin{thebibliography}{10}

\bibitem{MAsROB89}
Maurice Auslander and Ragnar-Olaf Buchweitz, \emph{The homological theory of
  maximal {C}ohen-{M}acaulay approximations}, M\'em. Soc. Math. France (N.S.)
  (1989), no.~38, 5--37, Colloque en l'honneur de Pierre Samuel (Orsay, 1987).
  \MR{MR1044344}

\bibitem{BCE20}
Silvana Bazzoni, Manuel Cort\'{e}s-Izurdiaga, and Sergio Estrada,
  \emph{Periodic modules and acyclic complexes}, Algebr. Represent. Theory
  \textbf{23} (2020), no.~5, 1861--1883. \MR{MR4140057}

\bibitem{BMPS19}
V\'{\i}ctor Becerril, Octavio Mendoza, Marco~A. P\'{e}rez, and Valente
  Santiago, \emph{Frobenius pairs in abelian categories. {C}orrespondences with
  cotorsion pairs, exact model categories, and {A}uslander-{B}uchweitz
  contexts}, J. Homotopy Relat. Struct. \textbf{14} (2019), no.~1, 1--50.
  \MR{MR3913970}

\bibitem{DLLM}
Zhenxing Di, Liping Li, Li~Liang, and Yajun Ma, \emph{Flat model structures and
  {G}orenstein objects in functor categories}, Proc. Roy. Soc. Edinburgh Sect.
  A \textbf{156} (2026), no.~2, 343--363. \MR{MR5047313}

\bibitem{rha}
Edgar~E. Enochs and Overtoun M.~G. Jenda, \emph{Relative homological algebra},
  de Gruyter Expositions in Mathematics, vol.~30, Walter de Gruyter \& Co.,
  Berlin, 2000. \MR{MR1753146}

\bibitem{Gi04}
James Gillespie, \emph{The flat model structure on {${\rm Ch}(R)$}}, Trans.
  Amer. Math. Soc. \textbf{356} (2004), no.~8, 3369--3390. \MR{MR2052954}

\bibitem{Gil162}
\bysame, \emph{Hereditary abelian model categories}, Bull. Lond. Math. Soc.
  \textbf{48} (2016), no.~6, 895--922. \MR{MR3608936}

\bibitem{RT12}
R\"{u}diger G\"{o}bel and Jan Trlifaj, \emph{Approximations and endomorphism
  algebras of modules. {V}olume 2}, extended ed., De Gruyter Expositions in
  Mathematics, vol.~41, Walter de Gruyter GmbH \& Co. KG, Berlin, 2012,
  Predictions. \MR{MR2985654}

\bibitem{HJ19}
Henrik Holm and Peter J{\o}rgensen, \emph{Cotorsion pairs in categories of
  quiver representations}, Kyoto J. Math. \textbf{59} (2019), no.~3, 575--606.
  \MR{MR3990178}

\bibitem{HJ19b}
\bysame, \emph{Model categories of quiver representations}, Adv. Math.
  \textbf{357} (2019), 106826, 46. \MR{MR4013804}

\bibitem{HJ21}
\bysame, \emph{The {$Q$}-shaped derived category of a ring}, J. Lond. Math.
  Soc. (2) \textbf{106} (2022), no.~4, 3263--3316. \MR{MR4524199}

\bibitem{HJ23}
\bysame, \emph{The {$Q$}-shaped derived category of a ring---compact and
  perfect objects}, Trans. Amer. Math. Soc. \textbf{377} (2024), no.~5,
  3095--3128. \MR{MR4744776}

\bibitem{Ho02}
Mark Hovey, \emph{Cotorsion pairs, model category structures, and
  representation theory}, Math. Z. \textbf{241} (2002), no.~3, 553--592.
  \MR{MR1938704}

\bibitem{LY22}
Li~Liang and Gang Yang, \emph{Constructions of {F}robenius pairs in abelian
  categories}, Mediterr. J. Math. \textbf{19} (2022), no.~2, Paper No. 76, 17.
  \MR{MR4384774}

\bibitem{Oberst70}
Ulrich Oberst and Helmut R\"{o}hrl, \emph{Flat and coherent functors}, J.
  Algebra \textbf{14} (1970), 91--105. \MR{MR257181}

\bibitem{Ru10}
Wolfgang Rump, \emph{Flat covers in abelian and in non-abelian categories},
  Adv. Math. \textbf{225} (2010), no.~3, 1589--1615. \MR{MR2673740}

\bibitem{SV26}
Anastasios Slaftsos and Jorge Vit\'{o}ria, \emph{Relative {Q}-shaped
  homological algebra}, preprint \textbf{\arxiv[RT]{2602.22986}}.

\bibitem{St68}
Bo~Stenstr\"{o}m, \emph{Purity in functor categories}, J. Algebra \textbf{8}
  (1968), 352--361. \MR{MR229697}

\bibitem{SS20}
Jan \v{S}aroch and Jan \v{S}t'ov\'{\i}\v{c}ek, \emph{Singular compactness and
  definability for {$\Sigma$}-cotorsion and {G}orenstein modules}, Selecta
  Math. (N.S.) \textbf{26} (2020), no.~2, Paper No. 23. \MR{MR4076700}

\end{thebibliography}

\def\cprime{$'$}
  \providecommand{\arxiv}[2][AC]{\mbox{\href{http://arxiv.org/abs/#2}{\sf
  arXiv:#2 [math.#1]}}}
  \providecommand{\oldarxiv}[2][AC]{\mbox{\href{http://arxiv.org/abs/math/#2}{\sf
  arXiv:math/#2
  [math.#1]}}}\providecommand{\MR}[1]{\mbox{\href{http://www.ams.org/mathscinet-getitem?mr=#1}{#1}}}
  \renewcommand{\MR}[1]{\mbox{\href{http://www.ams.org/mathscinet-getitem?mr=#1}{#1}}}
\providecommand{\bysame}{\leavevmode\hbox to3em{\hrulefill}\thinspace}
\providecommand{\MR}{\relax\ifhmode\unskip\space\fi MR }
\providecommand{\MRhref}[2]{%
  \href{http://www.ams.org/mathscinet-getitem?mr=#1}{#2}
}
\providecommand{\href}[2]{#2}

\end{document}